\documentclass[12pt,oneside]{amsart}
\usepackage{epsfig}
\usepackage{graphicx}
\usepackage{amscd,amssymb}
\theoremstyle{plain}
\newtheorem{Thm}[equation]{Theorem}

\newtheorem{Rmk}[equation]{Remark}
\newtheorem{Cor}[equation]{Corollary}
\newtheorem{Con}[equation]{Conjecture}
\newtheorem{Prop}[equation]{Proposition}
\newtheorem{Lem}[equation]{Lemma}

\numberwithin{equation}{section}

\newcommand{\mc}[1]{{}}

\newcommand{\e}{\varepsilon}

\newcommand{\br}{\mathbb{R}}

\newcommand{\Di}{\mathbb{D}}

\newcommand{\ba}{\backslash}

\newcommand{\G}{\Gamma}

\newcommand{\Cal}{\mathcal}

\newcommand{\bi}{\begin{itemize}}

\newcommand{\ei}{\end{itemize}}

\newcommand{\Vol}{\operatorname{Vol}}
\newcommand{\vol}{\operatorname{Vol}}

\newcommand{\T}{\operatorname{T}}
\newcommand{\SL}{\operatorname{SL}}

\newcommand{\Ad}{\operatorname{Ad}}

\newcommand{\vs}{\vskip 5pt}

\begin{document}

\dedicatory{Dedicated to Hillel Furstenberg with respect and admiration}
\title[orbits of discrete subgroups]{Orbits of discrete subgroups on a symmetric space and the Furstenberg boundary}
\author{Alexander Gorodnik and Hee Oh}
\address{Mathematics department\\ University of Michigan\\  Ann Arbor, MI 48109}
\email{gorodnik@umich.edu}

\address{Mathematics 253-37\\
Caltech\\Pasadena, CA 91106}
\email{heeoh@caltech.edu}

\thanks{The first author is partially supported by NSF grant 0400631.}
\thanks{The second author partially supported by NSF grant 0333397.}

 \abstract{Let $X$ be a symmetric
space of noncompact type and $\Gamma$ a lattice in
 the isometry group of $X$. We study the distribution of orbits of $\G$ acting on the
 symmetric space $X$ and its geometric boundary $X(\infty)$. More precisely, for
 any $y\in X$ and $b\in X(\infty)$, we investigate the distribution of the set
 $\{(y\gamma, b\gamma^{-1}):\gamma \in \G\}$ in $X\times X(\infty)$.
 It is proved, in particular, that the orbits of $\Gamma$ in the Furstenberg boundary
 are equidistributed,  and that the orbits of $\Gamma$ in $X$ are equidistributed in
 ``sectors'' defined with respect to a Cartan decomposition. We also discuss an application to
 the Patterson-Sullivan theory.
Our main tools are the strong wavefront lemma and the equidistribution of solvable flows
 on homogeneous spaces.
}\endabstract 

\maketitle

{\small \tableofcontents
}

\section{Introduction}

\begin{figure}[b] \label{pic1}
\begin{center}
\includegraphics[clip=true,width=7cm,height=7cm]{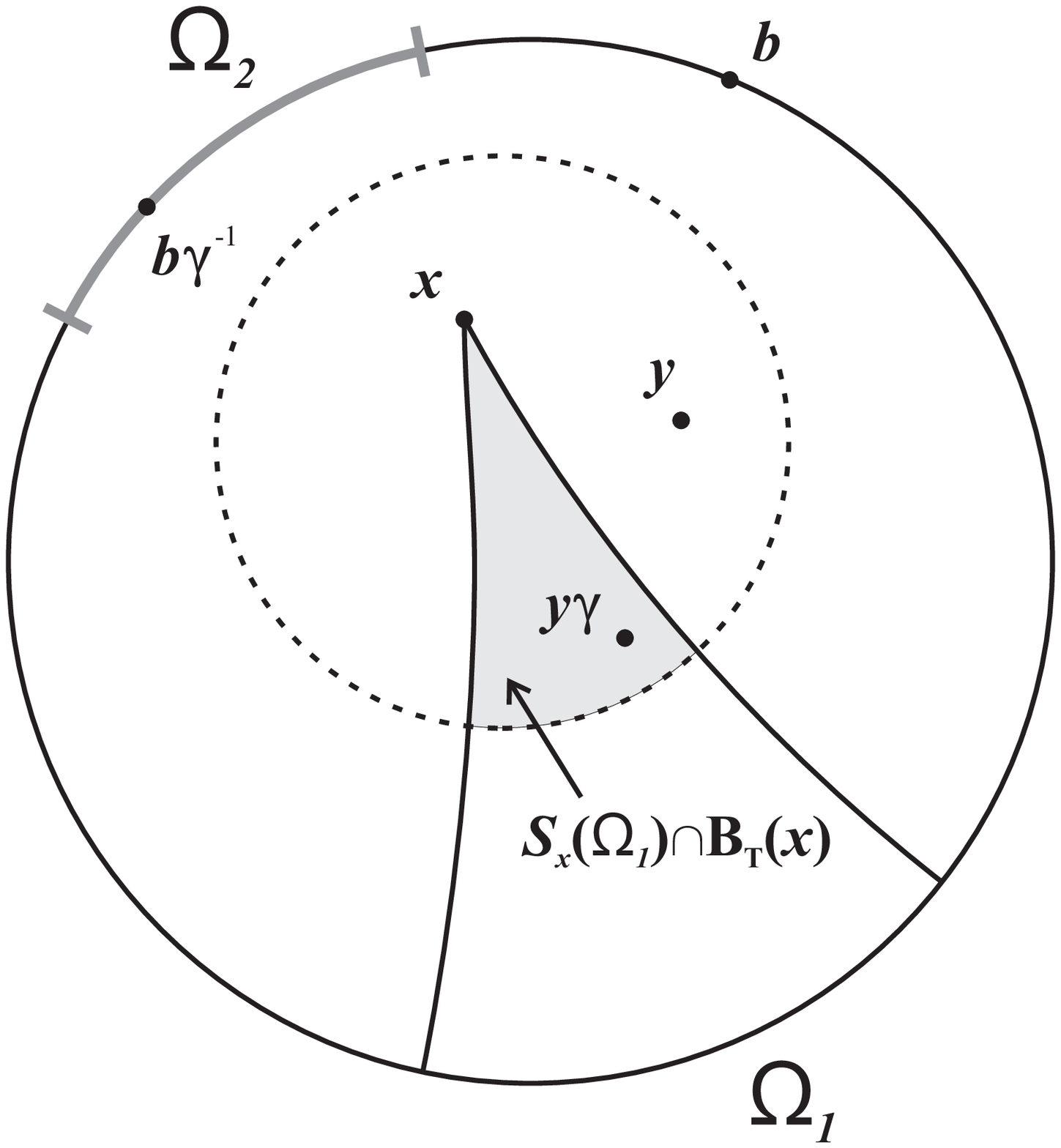}
\end{center}
\caption{}
\end{figure}

Let $\Di$ denote the hyperbolic unit disc and $\G$ a torsion free
discrete subgroup of the isometry group of $\Di$ such that
$\Di/\G$ has finite area. The geometric boundary of $\Di$ is the
space of the equivalence classes of geodesic rays in $\Di$. 
It can be identified with the unit circle $\mathbb{S}$.
Note that the action of $\Gamma$ on $\Di$ extends to the geometric boundary of
$\Di$.

Let $x\in \mathbb{D}$. We denote by $\operatorname{B_T}(x)$ the ball of radius $T$ centered at $x$. For an arc
$\Omega\subset \mathbb{S}$, the sector
$\Cal S_x(\Omega)$ in $\Di$ is defined to be the set of
points $z\in \mathbb{D}$ such that the end point of the geodesic
ray from $x$ to $z$ lies in $\Omega$. Denote by $m_x$ the unique
probability measure on $\mathbb{S}$ invariant under the isometries that
fix the point $x$. Then:

\begin{itemize}
\item[A.] For any $x,y\in \Di$, $b\in \mathbb{S}$ and an arc $\Omega\subset \mathbb{S}$,
$$
\#\{\gamma\in\Gamma:\, b\gamma^{-1}\in \Omega,\, y\gamma\in
\operatorname{B_T}(x) \}\sim_{T\to \infty}
m_y(\Omega)\cdot\frac{\hbox{\rm Area}(\operatorname{B_T}(x))}{\hbox{Area}(\Di/\G)} .$$
\item[B.]
For any $x,y\in\Di$ and an arc $\Omega\subset\mathbb{S}$,
$$
\#\{\gamma\in \G: y\gamma
\in \Cal S_x(\Omega)\cap \operatorname{B_T}(x)
\}\sim_{T\to \infty} m_x(\Omega)\cdot\frac{\hbox{\rm Area}
(\operatorname{B_T}(x))}
{\hbox{Area}(\Di/\G)}.
$$
\item[C.]
For every $x,y\in\Di$, $b\in \mathbb{S}$ and arcs $\Omega_1,\Omega_2\subset\mathbb{S}$,
\begin{multline*}
\#\{\gamma\in\Gamma:\, y\gamma\in \Cal S_x(\Omega_1)\cap\operatorname{B_T}(x),\, b\gamma^{-1}\in\Omega_2\}\\
\sim_{T\to \infty} m_x(\Omega_1)m_y(\Omega_2)\cdot\frac{\hbox{\rm Area} (\operatorname{B_T}(x))}{\hbox{Area}(\Di/\G)}.
\end{multline*}
\end{itemize}

Statement (A) may be deduced from the work of A. Good \cite{good}. 
Statement (B) was shown by P. Nicholls \cite{Ni} (see also \cite{S}).
Statement (C), which is new, shows that the
equidistribution phenomena in (A) and (B) are indeed independent.

The main purpose of this paper is to obtain an analog of statement (C)
(note that (C) implies both (A) and (B)) for an arbitrary
Riemannian symmetric space of noncompact type (see Theorems \ref{thm_main0} and \ref{thm_main} below).
We also generalize statement (B) to the equidistribution of lattice points in a connected noncompact
semisimple Lie group $G$ with finite center
with respect to both $K$-components in a Cartan decomposition $G=KA^+K$ (see Theorem \ref{th_KK} below).

\vs
Let $X$ be a Riemannian symmetric space of noncompact type
and $X(\infty)$ the geometric boundary of $X$ (that is, the space
of equivalence classes of geodesic rays in $X$). Denote by $G$ the
identity component of the isometry group of $X$ acting on $X$ from
the right. Let $\Gamma$ be a lattice in $G$, i.e., a discrete
subgroup with finite covolume. The action of $G$
on $X$ extends to $X(\infty)$.

For $x\in X$, we denote by $\operatorname{B_T}(x)$ the Riemannian
ball of radius $T$ centered at $x$, by $K_x$  the stabilizer of $x$ in $G$, and by $\nu_{x}$
the probability Haar measure on $K_x$. For $x\in X$ and $b\in X(\infty)$, we denote by
$m_{b,{x}}$ the unique probability ${K}_x$-invariant measure
supported on the orbit $bG\subset X(\infty)$ (note that $G$ acts
transitively on $X(\infty)$ only when the rank of $X$ is one, and
that $K_x$ acts transitively on each $G$-orbit in $X(\infty)$).
Fix a closed Weyl chamber $\mathcal{W}_{x} \subset X$ at $x$.
According to the Cartan decomposition, we have $X=\mathcal{W}_{x}
K_x$. Let $M_x$ denote the stabilizer of $\mathcal{W}_{x}$ in
$K_x$.

The following is one of our main theorems:
\begin{Thm}\label{thm_main0} For $x, y \in X$,
$b\in  X(\infty)$, and any Borel subsets  $\Omega_1\subset
K_x$ and $\Omega_2\subset bG$ with boundaries of measure zero,
\begin{equation*}
\#\{\gamma\in\Gamma:\, y\gamma\in \mathcal{W}_{x}\Omega_1\cap
\operatorname{B_T}(x),\, b\gamma^{-1}\in\Omega_2\}\sim_{T\to
\infty} \nu_{x}(M_x\Omega_1)m_{b, y}(\Omega_2)\cdot
\frac{\hbox{\rm Vol}(\operatorname{B_T})}{\Vol(G/\G)},
\end{equation*}
where $\hbox{\rm Vol}(\operatorname{B_T})$ denotes the volume of a ball of radius $T$ in $X$.
\end{Thm}

We deduce Theorem \ref{thm_main0} from a stronger result on the level of Lie groups.
Fix the following data:
\begin{itemize}
\item $G$ --  a connected noncompact
semisimple Lie group with finite center,
\item $G_n$ -- the product of all noncompact simple factors of $G$,
\item $G=K_1A^+K_1$ -- a Cartan decomposition of $G$,
\item $d$ -- an invariant metric on $K_1\ba G$,
\item $K_2$ -- a maximal compact subgroup of $G$,
\item $Q$ --  a closed subgroup of $G$
that contains a maximal connected split solvable subgroup.
\end{itemize}
Recall that a solvable subgroup $S$ is called {\it split} if the eigenvalues of any element of $\hbox{\rm Ad}(S)$ are real for the
adjoint representation $\hbox{\rm Ad}:G\to \hbox{GL}(\hbox{Lie}(G))$.
It is well-known that a maximal connected split solvable subgroup is a subgroup of the form $AN$ for an Iwasawa decomposition $G=KAN$.
Thus, $G={K}_2Q$. Denote by $\nu_1$ and $\nu_2$ the probability Haar measure on $K_1$ and $K_2$ respectively.
Let $M_1$ be the centralizer of $A$ in $K_1$ and $M_2=K_2\cap Q$.
Since any two maximal compact subgroups of $G$ are
conjugate to each other, there exists $g\in G$ such that $K_2=g^{-1} K_1 g$.
Let $\G$ be a lattice in $G$ such that $\overline{\Gamma G_n Q}=G$.

\begin{Thm}\label{thm_main}\label{main}
For any Borel subsets $\Omega_1\subset  K_1$ and $\Omega_2\subset K_2$
with boundaries of measure zero,
\begin{multline*}
\#\{\gamma\in\Gamma\cap g^{-1}K_1{A^+}\Omega_1\cap
\Omega_2 Q:\, d(K_1,K_1g\gamma)<T\}\\
\sim_{T\to\infty}
\nu_1(M_1\Omega_1)\nu_2(\Omega_2M_2)\cdot \frac{\hbox{\rm
Vol}(G_T)}{\Vol(G/\Gamma)},
\end{multline*}
where $\hbox{\rm Vol}(G_T)$
denotes the volume of a Riemannian ball of radius $T$ in $G$.
\end{Thm}

To understand the presence of $M_1$ and $M_2$ in the above
asymptotics, observe that $K_1A^+\Omega_1 =K_1A^+M_1\Omega_1$ and $\Omega_2M_2Q=\Omega_2Q$.

\begin{Rmk} {\rm We mention that the continuous version of Theorem \ref{main} does not seem obvious either.
The method of the proof of Theorem \ref{thm_main} also yields the following volume asymptotics:
\begin{multline*}
\Vol(\{h\in g^{-1}K_1{A^+}\Omega_1\cap \Omega_2  Q:\,
d(K_1,K_1gh)<T\})\\ \sim_{T\to\infty}
\nu_1(M_1\Omega_1)\nu_2(\Omega_2M_2)\cdot \frac{\hbox{\rm
Vol}(G_T)}{\Vol(G/\Gamma)}.
\end{multline*}}
\end{Rmk}

\begin{Rmk}\label{r_one_side}
{\rm
If in Theorem \ref{thm_main} we replace $K_1{A^+}\Omega_1$ by $\Omega_1{A^+}K_1$, then the statement of
the theorem is false.
In fact,
we can show that there exist nonempty open subsets $\Omega_1\subset K_1$ and
$\Omega_2\subset K_2$ such that
\begin{equation*}
\lim_{T\to\infty}\frac{1}{\Vol(G_T)} \Vol (\{h\in g^{-1}\Omega_1{A^+}K_1\cap \Omega_2
Q:\,d(K_1,K_1gh)<T\})= 0.
\end{equation*}
}
\end{Rmk}

To state yet another
generalization of statement (B), we fix a Cartan decomposition $G=KA^+K$
and an invariant Riemannian metric $d$ on $K\ba G$.
Let $\G$ be any lattice in $G$.
Recall that Eskin and McMullen showed in \cite{EM} that
for a lattice $\G$ in $G$,
\begin{equation}\label{eq_em}
\#\{\gamma\in\G:d(K,Kg\gamma)<T\}=\# (\G\cap g^{-1}KA^+_T K)\sim_{T\to \infty} \frac{\Vol (G_T)}{\Vol(G/\G)},
\end{equation}
where $A^+_T=\{a\in A^+: d(K, Ka)< T\}$.
The following theorem generalizes this result:
\begin{Thm}\label{th_KK} For $g\in G$ and any Borel subsets $\Omega_1\subset K $ and $\Omega_2\subset K$
with boundaries of measure zero,
\begin{eqnarray*}
\#(\Gamma\cap g^{-1}\Omega_1A_T^+M\Omega_2)\sim_{T\to\infty} \frac{\Vol(g^{-1}\Omega_1A_T^+M\Omega_2)}{\Vol
(G/\Gamma)}
=\nu(\Omega_1M)\nu(M\Omega_2)\cdot \frac{\Vol(G_T)}{\Vol
(G/\Gamma)},
\end{eqnarray*}
where $M$ is the centralizer of $A^+$ in $K$, and $\nu$ is the
probability Haar measure on $K$.
\end{Thm}

We now present several corollaries of (the methods of) the
above theorems. 

\subsection{Lattice action on the Furstenberg boundary}
For a connected semisimple Lie group $G$ with finite center, the
Furstenberg boundary of $G$ is identified with the quotient space
$G/P$ where $P$ is a minimal parabolic subgroup of $G$ (see
\cite[Ch.~IV]{gjt}). In the rank one case, the Furstenberg
boundary $G/P$ coincides with the geometric boundary $X(\infty)$
of the symmetric space $X$ of $G$. In the higher rank case, $G/P$ is
isomorphic to the $G$-orbit in $X(\infty)$ of any regular geodesic
class and can be identified with the space of asymptotic
classes of Weyl chambers in $X$.

It is well-known
that the action of a lattice $\G$ on $G/P$ is minimal, i.e., every $\G$-orbit
is dense (\cite[Lemma~8.5]{Mo}). A natural question is whether each $\G$-orbit
in $G/P$ is equidistributed.
Corollary \ref{P}, which is a special case of Theorem \ref{main}, implies an affirmative answer
in a much more general setting.

Let $d$ denote an invariant Riemannian metric on the symmetric space $X\simeq K\ba G$,
where $K$ is a maximal compact subgroup of $G$.

\begin{Cor}\label{P} Let $Q$ be a closed subgroup of $G$ containing a maximal connected split solvable subgroup of $G$
and $g\in G$. Denote by $\nu_g$ the
unique $g^{-1} Kg$-invariant probability measure on $G/Q$.
Let $b\in G/Q$ and $\G$ be a lattice in $G$ such that $\overline{\Gamma G_n b}=G/Q$.
Then for any Borel subset
$\Omega\subset G/Q$ such that $\nu_g(\partial\Omega)=0$,
$$\#\{\gamma\in \Gamma:\,\gamma b\in \Omega ,\,
d(K,Kg\gamma)<T\}\sim_{T\to \infty} \nu_g(\Omega)\cdot \frac{\vol
(G_T)}{\vol (G/\G)}.$$
\end{Cor}

It follows from the result of N.~Shah \cite[Theorem~1.1]{Sh} and Ratner's topological rigidity
\cite{Ra2} that the condition $\overline{\Gamma G_n b}=G/Q$ is equivalent to the density of the orbit $\Gamma b$ in $G/Q$.



In the case when $Q$ is a parabolic subgroup of $G$, a different proof of this result, which is based on ideas
developed in \cite{mau}, is recently given in \cite{gm}.

In the last decade or so there have been intensive studies on
the equidistribution properties of lattice points on homogeneous
spaces of $G$ using various methods from analytic number theory,
harmonic analysis and ergodic theory (\cite{DRS}, \cite{EM},
\cite{EMM},
\cite{EMS}, \cite{GO}, \cite{EO}, \cite{Go}, \cite{l}, \cite{mau},
\cite{no}, etc.). Of particular interest is the case when the
homogeneous space is a real algebraic variety. While most of the attention
in this direction is focused on the case of affine homogeneous
varieties, there is not so much work done for the projective
homogeneous varieties, except for recent works \cite{Go} and
\cite{mau}. In \cite{Go}, one studies the distribution of lattice
orbits on the real projective homogeneous varieties of
$G=\SL_n(\br)$ with respect to the norm given by $\|g\|=\sqrt{\sum
g_{ij}^2}$, $g\in\SL_n(\br)$. In \cite{mau}, one investigates the
distribution of lattice orbits on the boundary of a real hyperbolic
space. Corollary \ref{P} extends both results
by proving that an orbit of a lattice in a connected
noncompact semisimple real algebraic group $G$ is equidistributed on
any projective homogeneous variety of $G$ (with respect to a
Riemannian metric). 

More generally, we state the following conjecture:

\begin{Con}
Let $G$ be a connected semisimple Lie group with finite center, $\Gamma$ a lattice in $G$,
and $Y$ a compact homogeneous space of $G$. Then every dense orbit of $\Gamma$ in $Y$ is equidistributed, i.e.,
there exists a smooth measure $\nu$ on $Y$ such that for any $y\in Y$ 
with $\overline{\Gamma y}=Y$ and for any Borel set $\Omega\subset Y$ with boundary of measure zero,
$$
\#\{\gamma\in \Gamma:\,\gamma y\in \Omega ,\,
d(K,K\gamma)<T\}\sim_{T\to \infty} \nu(\Omega)\cdot \frac{\vol
(G_T)}{\vol (G/\G)}.
$$
\end{Con}

The structure of compact homogeneous spaces of $G$ was studied in \cite{wit}.
We note that the case of the conjecture when $Y=G/\G$ for a cocompact lattice is also known (see Theorem \ref{comme}).

\subsection{Measure-preserving lattice actions}
\label{sec_comme}
Let $G$ be a connected semisimple noncompact Lie group with finite center and $\Gamma_1$, $\Gamma_2$
lattices in $G$. We consider the action of $\Gamma_1$ on $G/\Gamma_2$.
Let $d$ be an invariant Riemannian metric on the symmetric space $K\ba G$.

\begin{Thm}\label{comme} 
Suppose that for $y\in G/\Gamma_2$, the orbit $\Gamma_1y$ is dense in $G/\Gamma_2$
Then for any $g\in  G$ and  any Borel subset $\Omega\subset G/\G_2$ with boundary of measure zero,
$$\{\gamma\in \G_1 : \gamma y \in \Omega, \, d(K, Kg \gamma)<T\}
\sim_{T\to \infty} \frac{\Vol (\Omega)\cdot {\vol(G_T)}}{\vol(
G/\G_1)\vol ( G/\G_2)}, $$ where all volumes are computed with
respect to one fixed Haar measure on $G$.
\end{Thm}

For example, Theorem \ref{comme} applies to the case when $G$ is a simple connected noncompact Lie group,
and $\Gamma_1$ and $\Gamma_2$ are noncommensurable lattices in $G$.
(Recall that the lattices are called {\it commensurable} if $\Gamma_1\cap \Gamma_2$ has finite index in both $\Gamma_1$ and $\Gamma_2$.)
It was first observed by Vatsal (see \cite{Va}) that $\G_1 \G_2$ is dense in $G/\G_2 $.
This is a (simple) consequence of Ratner's topological rigitity \cite{Ra2}. 

Theorem \ref{comme} was
proved in \cite{Oh} for $G=\SL_n(\br)$ equipped with the norm
$\|g\|=\sqrt{\sum g_{ij}^2}$ and in \cite{GW} for general semisimple
Lie groups without compact factors.

\subsection{Counting lattice points lying in sectors and an
application to the Patterson-Sullivan theory} Keeping the
notations of Theorem \ref{thm_main0}, we note that $\Cal W_{x}
\Omega_1$ is an analog of the sector $\Cal
S_x(\Omega_1)$ discussed in the beginning of the introduction.
Hence the following corollary, which is a special case of Theorem \ref{th_KK}, is a generalization of (B):
\begin{Cor}\label{limit} For $x, y\in X$ and a Borel subset $\Omega\subset K_x$
with boundary of measure zero,
$$ \# (y \G \cap \Cal W_{x} \Omega \cap \operatorname{B_T}(x) )
 \sim_{T\to \infty} \#(\Gamma\cap K_y)\cdot
\nu_{x}(M_x\Omega)\cdot\frac{\vol (\operatorname{B_T})}{\vol
(G/\G)} .$$
\end{Cor}
This result was shown by Nicholls \cite{Ni} for the case of
a real hyperbolic space (see also \cite{S}).

 Corollary \ref{limit} is related to the theory of
Patterson-Sullivan measures. These measures, which were introduced
by Patterson \cite{pat} and Sullivan \cite{sul} in the case of a
real hyperbolic space, have proved to be an invaluable tool for
the study of spacial distribution of orbits of discrete groups of isometries.
The theory of Patterson and Sullivan was extended to higher rank
symmetric spaces by several authors (see \cite{Al}, \cite{qui},
\cite{link}).

Recall that $X\cup X(\infty)$ equipped with conic topology is one
of the standard compactifications of $X$ (see \cite[Ch.~III]{gjt}).
Denote by $\delta_\Gamma$ the critical exponent of the Dirichlet series $\sum_{\gamma\in\Gamma}e^{-sd(x,y\gamma)}$ where
$x,y\in X$. For $x,y\in X$ and $s>\delta_\Gamma$, define measures on $X\cup X(\infty)$:
\begin{equation}\label{eq_mxy}
\mu_{x,y,s}=\frac{1}{\sum_{\gamma\in\Gamma}e^{-sd(x,y\gamma)}}\sum_{\gamma\in\Gamma}e^{-sd(x,y\gamma)}D_{y\gamma},
\end{equation}
where $D_z$ denotes the Dirac measure at a point $z\in X$.
The family of Patterson-Sullivan measures $\{\mu_{x}: x\in X\}$ is
defined by 
$$
\mu_x=\lim_{s\to\delta_{\Gamma}^+} \mu_{x,y,s}.
$$
Using Corollary \ref{limit}, we deduce the following characterization of the Patterson-Sullivan measures:
\begin{Cor}\label{Albu}
For every $x\in X$, $\mu_{x}$ is the unique $K_x$-invariant probability measure supported on the set $bK_x\subset X(\infty)$,
where $b$ is the equivalence class containing the geodesic ray emanating from $x$ in the direction of the barycenter of the Weyl chamber $\mathcal{W}_x$.
\end{Cor}


Corollary \ref{Albu} was obtained by Albuquerque \cite{Al}, extending the
work of Patterson \cite{pat}, Sullivan \cite{sul} in the rank one case and
the work of Burger \cite{bu} in the case of the product of rank one spaces. Our
proof has a completely different flavor compared to their geometric
methods.

\subsection{Acknowledgments} The first author wishes to thank
for the hospitality of the Department of Mathematics in California
Institute of Technology, where most of this work has been done.

\section{Main ingredients of the proofs}
\subsection{The strong wavefront lemma}  The following theorem is a
basic tool which enables us to reduce the counting problems for
$\Gamma$ as in Theorems \ref{main} and \ref{th_KK} to the study of
continuous flows on the homogeneous space $\Gamma\ba G$.

Let $G$ be a connected noncompact semisimple Lie group with finite center,
$G=KA^+K$ a Cartan decomposition, and $M$ the centralizer of $A$ in $K$.

\begin{Thm}[\bf The strong wavefront lemma] \label{wave}
Let $\Cal C$ be any closed subset of $A^+$ with a positive distance from the walls of $A^+$.
Then for any neighborhoods $U_1, U_2$ of $e$ in $K$ and $V$
of $e$ in $A$, there exists a neighborhood
 $\Cal O $ of $e$ in $G$
such that for any $g=k_1ak_2\in K \Cal C K$,
 \begin{itemize}
 \item[(1)]
 $g \Cal O \subset
(k_1U_1)(a VM)(k_2U_2)$;
\item[(2)]
 $\Cal O g \subset
(k_1U_1)(a VM)(k_2U_2) $. \end{itemize}
\end{Thm}

\begin{Rmk} {\rm
One can check that Theorem \ref{wave} fails if the set $\mathcal{C}$ contains a sequence that converges to a point
in a wall of the Weyl chamber $A^+$.
}
\end{Rmk}

Theorem \ref{wave} has several geometric implications for the
symmetric space $K\ba G$:
\begin{itemize}
\item {\it Strengthening of the wavefront lemma.} Recall that
the wavefront lemma introduced by Eskin and McMullen in \cite{EM} says
that for any neighborhood $\Cal O'$ of $e$ in $G$, there exists a
neighborhood $\Cal O$ of $e$ in $G$ such that
\begin{equation}\label{em} a\Cal O  \subset  \Cal O' a K\quad\hbox{for all $a\in A^+$}.
\end{equation}
To see that our strong wavefront lemma (1) implies the wavefront lemma for $a\in A^+$ 
with at least a fixed positive distance from the walls of $A^+$, note that $\Cal O'$ contains $U_1V$ for some
neighborhood $U_1$ of $e$ in $K$ and some neighborhood $V$ of $e$
in $A$, and hence
$$
U_1 a
VMK = U_1 V a K\subset  \Cal O' a K
$$
By Theorem \ref{wave} (1), there exists a
neighborhood $\Cal O$ of $e$ such that
$$ a \Cal O\subset U_1 a VMK.$$
Thus, $a \Cal O\subset \Cal O' aK$.

\begin{figure}[t] \label{pic2}
\begin{center}
\includegraphics[clip=true,width=7cm,height=7cm]{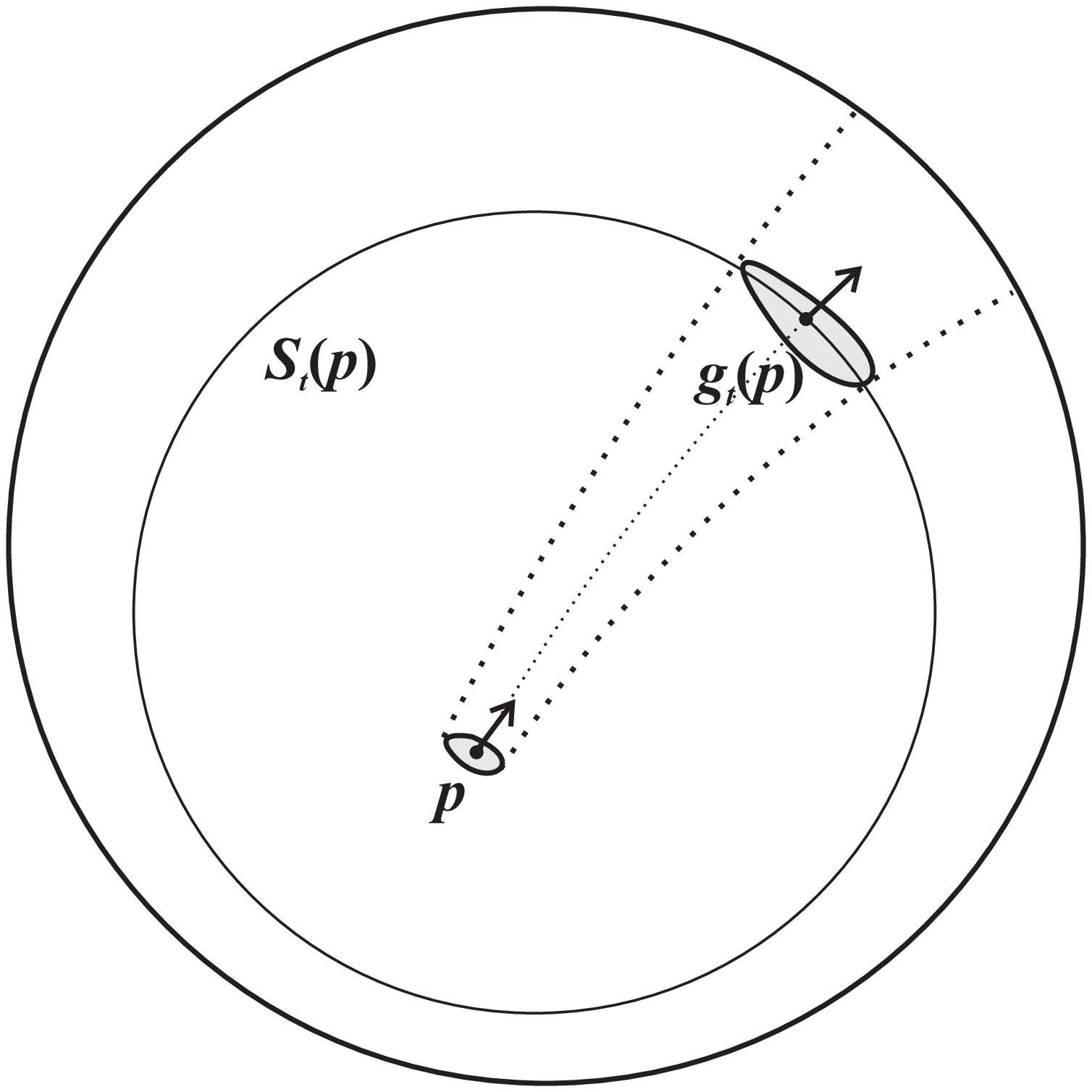}
\end{center}
\caption{}
\end{figure}

To illustrate geometric meaning of the strong wavefront lemma (1), we consider the
unit disc $\Di$ equipped with the standard hyperbolic metric and the geodesic flow $g_t$ on the
unit tangent bundle $\T^1(\Di)$ that transports a vector distance $t$ along the
geodesic to which it is tangent. Note that with the identification $\T^1(\Di)\simeq \hbox{PSL}_2(\mathbb{R})$,
the geodesic flow $g_t$ corresponds to the left multiplication by $\left(
\begin{tabular}{cc}
$e^{t/2}$ & $0$\\
$0$ & $e^{-t/2}$
\end{tabular}
\right)$.
Let $p\in \Di$ and let $K\subset
\T^1(\Di)$ be the preimage of $p$ under the projection map $\pi:\T^1(\Di)
\to \Di$. Note that $K$ consists of vectors lying over $p$ and
pointing in all possible directions, and $g_t(K)$ consists of the unit
vectors normal to the sphere $S_t(p)\subset \Di$ of radius $t$.
The wavefront lemma (see (\ref{em})) implies that one can find a
neighborhood $\Cal O\subset \T^1(\Di)$ of a vector $v$ based at $p$
such that $g_t(\Cal O)$ remains close to $g_t(K)$ uniformly for every $t\ge 0$.

However, this does not compare $g_t(\Cal O)$ with the vector
$g_t(v)$, but rather with the set $g_t(K)$. Theorem \ref{wave} (1) says
that we may choose a neighborhood $\Cal O$ of $v$ in $\T^1(\Di )$
so that $g_t(\Cal O)$ is close to the vector $g_t(v)$ uniformly on
$t$ in both angular and radial components (see Fig. 2).

\vs\item {\it Uniform openness of the map $K\times A^+ \times K\to G$.}
 The product map
$$K\times (\hbox{interior of }A^+) \times K \to G$$ is  a diffeomorphism onto
a dense open subset in $G$, and in particular, it is an open map.
Theorem \ref{wave} (2) shows that this map is {\it uniformly open}
with respect to the base of neighborhoods $\mathcal{O}g$, where 
$\mathcal{O}$ is a neighborhood of $e$ in $G$ and $g\in G$,
on any subset contained in $K\times A^+ \times K$ with a positive distance from the walls of $A^+$.

We illustrate the geometric meaning of this property for the case of the hyperbolic unit disc $\Di$.
It follows from Theorem \ref{wave} (2) that for every neighborhood $\mathcal{O}'$ of $e$ in $G$,
there exists a neighborhood $\mathcal{O}$ of $e$ in $G$ such that for every $a\in A^+$ with at least a fixed positive
distance from the walls of $A^+$,
$$
\mathcal{O}a\subset Ka\mathcal{O}'.
$$
This implies that for every open subset $\mathcal{O}'\subset \T^1(\mathbb{D})$, $p\in\mathcal{O}'$, and any $t>0$, the set
$\pi(g_t(\mathcal{O}'))$ contains a ball centered at $\pi(g_t(p))$ with a radius independent of $t$.

\vs\item {\it Well-roundness of bisectors.} Theorem
\ref{wave} (2) also implies the following corollary, which is a generalization of
the well-known property that the Riemannian balls 
 are well-rounded (a terminology used in \cite{EM}). Note that
the Riemannian ball $\{g\in G: d(K,Kg)<T\}$ is of the form $KA_T^+K$,
where $A_T^+=\{a\in A^+: d(K, Ka)< T\}$.

\begin{Cor}  For Borel subsets $\Omega_1, \Omega_2\subset K$ whose
boundary has measure zero, the family $ \{\Omega_1 A_T^+ \Omega_2:
T >0\} $ of bisectors is {\rm well-rounded}, i.e., for
every $\e>0$, there exists a neighborhood $\Cal O$ of $e$ in $G$
 such that
$$\Vol (\Cal O\cdot \partial (\Omega_1 A_T^+\Omega_2))\le
\e \cdot \Vol(\Omega_1 A_T^+ \Omega _2) $$ for all $T>0$.
\end{Cor}

To see how Theorem \ref{wave} (2) implies the above corollary, it
suffices to note that by the integration formula for a Cartan
decomposition (see (\ref{cartan}) below), there exist a small neighborhood
$U$ of $e$ in $K$ and a small neighborhood $V$ of $e$ in $A$
such that
$$\Vol ( (\partial (\Omega_1) U)( \partial ( A_T^+) VM)( \partial
(\Omega_2) U) ) \le \e \cdot \Vol (\Omega_1 A_T^+ \Omega_2 ).$$
\end{itemize}

\vs\subsection{Uniform distribution of
solvable flows} Using the strong wavefront lemma,
Theorem \ref{main} is deduced from
Theorem \ref{Pintro} below, which is also of
independent interest from the viewpoint of ergodic theory.

Let $K$ be a maximal compact subgroup of $G$ with the probability Haar measure $\nu$.
Let $Q$
be a closed subgroup of $G$ containing a maximal connected split solvable subgroup of $G$,
and let $\rho$ be a right invariant Haar measure on $Q$.
Fix a Cartan decomposition $G=KA^+K$ and $g\in G$. For $T>0$ and a subset $\Omega\subset K$, we define
$$Q_T(g, \Omega)= \{q\in Q: q\in  g^{-1}K{A^+}\Omega, \,
d(K, Kgq)< T\}.$$
If $\Omega=K$ and $g=e$, the set $Q_T(g, \Omega)$ is simply $\{q\in Q: d(K, Kq)< T\}$.
Recall that $G_n$ denotes the product of all noncompact simple factors of $G$.

\begin{Thm}\label{Pintro} 
Let $G$ be realized as a closed
subgroup of a Lie group $L$. Let $\Lambda$ be a lattice in $L$.
Suppose that for $y\in\Lambda\ba L$, the orbit $y G_n$ is dense in
$\Lambda\ba L$. Then for any Borel subset $\Omega\subset K$ with boundary of measure zero
and $f\in C_c(\Lambda\ba L)$,
\begin{equation*}\lim_{T\to \infty} \frac{1}{\rho (Q_T(g, K))}\int_{Q_T(g,
\Omega)} f(y q^{-1})\, d\rho(q)=
\frac{\nu(M\Omega)}{\mu(\Lambda\ba L)} \int_{\Lambda \ba L}f
\, d\mu,
\end{equation*}
where $M$ is the centralizer of $A$ in $K$, and $\mu$ is an $L$-invariant measure on $\Lambda\ba L$.
\end{Thm}

A main ingredient of the proof of Theorem \ref{Pintro} is the work
of N.~Shah \cite{Sh} (see Theorem \ref{shah} below) on the distribution in
$\Lambda\ba L$ of translates $ y U g$ as $g\to\infty$ for a subset $U\subset K$.
Shah's result is based on Ratner's
classification of measures invariant under unipotent flows \cite{Ra1}
and the work of Dani and Margulis on behavior of unipotent
flows \cite{DM}. Implementation of Shah's theorem in our setting is based on the
fundamental property of the Furstenberg boundary $\mathcal{B}$ of $G$: every
regular element in a positive Weyl chamber acts on an open
subset of full measure in $\mathcal{B}$ as a contraction.

\begin{Rmk}[\bf on the rate of convergence] {\rm The method
of the proof of Theorem \ref{shah} in \cite{Sh} does not give any
estimate on the rate of convergence. In the case when $L=G$ and
$U=K$, Theorem \ref{shah} was proved by Eskin and McMullen in
\cite{EM}. The latter proof is based on the decay of the matrix
coefficients of the quasi-regular representation of $G$ on
$L^2(\Gamma\ba G)$ and provides an estimate on the rate of
convergence. Combining the strong wavefront lemma (Theorem
\ref{wave}) with the method from \cite{EM}, we can derive an
estimate for the rate of convergence in Theorem \ref{Pintro} when
$L=G$ provided that one knows the rate of decay of matrix
coefficients of $L^2(\Gamma\ba G)$. In this case, it is also
possible to obtain rates of convergence for the theorems stated in
the introduction. We hope to address this problem in a sequel
paper.
}
\end{Rmk}

\vs
\subsection{Equidistribution of lattice points in bisectors} For
$\Omega_1, \Omega_2\subset K$, and $g\in G$, we define
$$G_T(g, \Omega_1,
\Omega_2)=\{h\in G: h\in g^{-1} \Omega_1 A^+ \Omega_2 ,\, d(K,
Kgh)<T\} .$$ Using the strong wavefront lemma (Theorem
\ref{wave} (2)), Theorem \ref{th_KK} is reduced to showing that the
sets $G_T(g,\Omega_1, \Omega_2)$ are equidistributed in
$\G\ba G$ in the sense of Theorem \ref{Pintro} for any Borel subsets $\Omega_1, \Omega_2\subset K$
with boundaries of measure zero.

\section{Cartan decomposition and the strong wavefront lemma}\label{strong2}

Let $G$ be a connected noncompact semisimple Lie group with finite center,
$K$ a maximal compact subgroup of $G$,
and $G=K\exp(\frak p)$ the Cartan decomposition determined by $K$.
A split Cartan subgroup $A$ with respect to $K$ is a maximal
connected abelian subgroup of $G$ contained in $\exp (\frak p)$. It
is well-known that two split Cartan subgroups with respect to $K$
are conjugate to each other by an element of $K$. Fix a split
Cartan subgroup $A$ of $G$ (with respect to $K$) with the set of
positive roots $\Phi^+$ and the positive Weyl chamber
$$A^+=\{a\in A: \alpha( \log a)\ge 0\text{ for all $\alpha\in
\Phi^+$}\} .$$
Set $\mathfrak{a}=\log (A)$ and $\mathfrak{a}^+=\log (A^+)$.
Let $M$ be the centralizer of $A$ in $K$.
Note that $M$ is finite if and only if $G$ is real split.

The following lemma is well-known (see, for example, \cite[Ch.~V]{Kn}).
\begin{Lem}\label{unique}
{\bf (Cartan decomposition)} For every $g\in G$, there exist a
unique element $\mu (g)\in \log(A^+)$ such that $g\in K \exp({\mu(g)})
K $. Moreover, if $k_1 ak_2= k_1'ak_2'$ for some $a$ in the
interior of $A^+$, then there exists $m\in M$ such that
$k_1=k_1'm$, $k_2=m^{-1} k_2'$.
\end{Lem}

Denote by $d$ an
invariant Riemannian metric
 on the symmetric space $K\backslash G$.

\begin{Lem} \label{l_hel}
For every $a_1$ and $a_2$ in the interior of ${A}^+$ and $k\in K$,
$$
d(Ka_1,Ka_2)\le d(Ka_1k,Ka_2).
$$
\end{Lem}

\begin{proof}
Let $a_i=\exp(H_i)$ for $H_i\in \mathfrak{a}^+$, $i=1,2$. Then $Ka_1k=K\exp(\hbox{Ad}(k^{-1})H_1)$.
Applying the cosine inequality (see \cite[Corollary~I.13.2]{hel})
 to the geodesic triangle with vertices
$Ke$, $Ka_1k$, and $Ka_2$, we obtain
\begin{eqnarray*}
d(Ka_1k,Ka_2)^2&\ge& d(K,Ka_1k)^2+d(K,Ka_2)^2-2d(K,Ka_1k)d(K,Ka_2)\cos\alpha\\
&=& \|H_1\|^2+\|H_2\|^2-2\|H_1\|\|H_2\|\cos\alpha,
\end{eqnarray*}
where $\alpha$ is the angle at the vertex $Ke$. Since
$$
\cos\alpha=\frac{\langle \Ad (k^{-1})H_1,H_2\rangle}{\|H_1\|\cdot\|H_2\|}
$$
and by Lemma \ref{link} below,
$$
\langle \Ad (k^{-1})H_1,H_2\rangle \le\langle H_1,H_2\rangle,
$$
it follows that
$$
d(Ka_1k,Ka_2)^2\ge \|H_1\|^2+\|H_2\|^2-2\langle  H_1,H_2\rangle
=\|H_1-H_2\|^2=d(Ka_1,Ka_2)^2.
$$
The lemma is proved.
\end{proof}

\begin{Lem}\label{link} For any $H_1$ and $H_2$ in the interior of $\frak a^+$ and
for any $k\in K$,
$$
\langle H_1,H_2\rangle \ge \langle \Ad (k)H_1,H_2\rangle .
$$
\end{Lem}
\begin{proof}
By \cite[Proposition~VIII.5.2]{hel} and its proof, every
$G$-invariant positive definite form on $K\ba G$ is of the form
$\sum_i \alpha_i B_i$, where $B_i$'s are the Killing forms of the
simple factors of $G$ and $\alpha_i >0$. Thus, it is sufficient to
consider the case when $G$ is simple, and the Riemannian metric is
given by the Killing form $B$ of $G$.

Define the function $f(k)=
 \langle \Ad(k)H_1,H_2\rangle $ on $K$.
Let $k_0\in K$ be a point where $f$ attains its maximum.
For every $Z\in\hbox{Lie}(K)$,
 \begin{align*}
0&=\frac{d}{dt}\Big|_{t=0} f(k_0 e^{tZ})=\frac{d}{dt}\Big|_{t=0} B(\Ad(k_0) \Ad(e^{tZ}) H_1, H_2) \\
&= B\left(\Ad(k_0)  \frac{d}{dt}\Big|_{t=0}( \Ad(e^{tZ}) H_1), H_2\right) =  B(\Ad(k_0) ( \hbox{ad}(Z ) H_1), H_2)\\
&=  B(\Ad(k_0) [Z, H_1], H_2) = B(\Ad(k_0)Z, [\Ad(k_0)H_1, H_2]). \end{align*}
This shows that $[\Ad(k_0)H_1, H_2]\perp \hbox{Lie}(K)$.
Since $[\Ad(k_0)H_1, H_2]\in \hbox{Lie}(K)$ and the restriction of
$B$ to $\hbox{Lie}(K)$ is negative definite, it follows that $[\Ad(k_0)H_1, H_2]=0$.
Therefore, $\Ad(k_0)H_1\in \frak a$.
Since the Weyl group $W$ acts transitively on the set of Weyl chambers
in $A$ and $K$ contains all representatives of the Weyl group,
there exists an element $w\in K$ that normalizes $A$ such that
$\Ad (w^{-1}k_0)H_1\in\mathfrak{a}^+$.
Since $H_1$ is in the interior of $\frak a^+$, it
follows from the uniqueness of the Cartan decomposition (Lemma \ref{unique})
that $\Ad (k_0)H_1=\Ad(w)H_1$.
It is easy to see \cite[p.~288]{hel} that $\|\Ad(w)H_1-H_2\|$, $w\in W$, achieves its
minimum at $w=e$. This implies that $\left<\Ad(w)H_1,H_2\right>$ is maximal for $w=e$
and finishes the proof.
\end{proof}

\begin{Prop}\label{strong}
Let $\Cal C$ be a closed subset contained in $A^+$ with a positive distance from the walls of $A^+$.
Then for any neighborhood $U_0$ of $e$ in $K$, there exists
$\varepsilon>0$ such that for any $a\in \Cal C$, $$ \{k\in K:
d(Kak,Ka)<\varepsilon\}\subset  MU_0.$$
\end{Prop}

\begin{proof}
Denote by $\Pi$ the set of simple roots corresponding to the Weyl chamber $A^+$.
Without loss of generality, we may assume that
$$
\mathcal{C}=\{a\in A^+: \alpha(\log a)\ge C\;\hbox{ for all $\alpha\in\Pi$}\}
$$
for some $C>0$. Suppose that
in contrary there exist sequences
$\{a_i\}\subset \Cal C$ and $\{k_i\}\subset K$ such that
$d(Ka_ik_i,Ka_i)\to 0$ as $i\to \infty$, and no limit points of
$\{k_i\}$ are contained in $M$. Passing to a subsequence,
we may assume that $k_i\to k_0$ as $i\to\infty$ for some
$k_0\in K-M$, and 
for every $\alpha\in\Phi$, the sequence $\{\alpha(\log a_i)\}$ is either bounded or divergent. Set
\begin{eqnarray*}
\Phi_{\pm}&=&\{\alpha\in\Phi:\, \alpha(\log a_i)\to\pm\infty\;\; \hbox{ as $i\to\infty$}\},\\
\Phi_{0}&=&\{\alpha\in\Phi:\, \{\alpha(\log a_i)\} \hbox{ is bounded}\}.
\end{eqnarray*}
Let $P^+$ be the standard parabolic subgroup
associated to $\Pi-\Phi_+$ and $P^-$ the standard opposite parabolic subgroup for $P^+$. Note that
$$
P^-=\{g\in G:\, \{a_iga_i^{-1}\} \hbox{ is bounded}\}.
$$
Denote by $U^+$ and $U^-$
the unipotent radicals of $P^+$ and $P^-$ respectively. Set $Z=P^+\cap P^-$,
so that $P^\pm=ZU^\pm$. It is easy to see that $P^\pm\cap K\subset Z$.
Denote by $\mathfrak{g}$, $\mathfrak{u}^+$, $\mathfrak{u}^-$, $\mathfrak{z}$ the corresponding Lie algebras
and by $\mathfrak{g}_\alpha$, $\alpha\in\Phi$, the root subspaces in $\mathfrak{g}$.
We have
\begin{equation}\label{eq_uzu}
\mathfrak{u}^\pm=\oplus_{\alpha\in\Phi_\pm}\, \mathfrak{g}_\alpha\quad\hbox{and}\quad
\mathfrak{z}=\oplus_{\alpha\in\Phi_{0}}\, \mathfrak{g}_\alpha.
\end{equation}

\vs\noindent\underline{Step 1}: We claim that $k_0\in Z$.

\vs There is an embedding $\pi: G/C_G\to\hbox{SL}_d(\mathbb{R})$, where $C_G$ is the center of $G$,
such that $\pi(AC_G)$ is contained in the group of diagonal matrices.
We have
$$
\pi(a_ik_ia_i^{-1})_{st}= \pi(a_i)_{ss}\pi(a_i)_{tt}^{-1}\cdot \pi(k_i)_{st},\quad s,t=1,\ldots,d.
$$
Passing to a subsequence, we may assume that for each $(s,t)$, the sequence $\{\pi(a_i)_{ss}\pi(a_i)_{tt}^{-1}\}$ is either bounded or
divergent. If the sequence $\{\pi(a_i)_{ss}\pi(a_i)_{tt}^{-1}\}$ is divergent, then $\pi(k_i)_{st}\to 0$ as $i\to\infty$. Thus, $\pi(k_0)_{st}=0$
for every pair $(s,t)$ such that $\{\pi(a_i)_{ss}\pi(a_i)_{tt}^{-1}\}$ is divergent. It follows that
$\pi(a_ik_0a_i^{-1})$ is bounded. Since the center of $G$ is finite, this proves that the sequence $\{a_ik_0a_i^{-1}\}$ is bounded.
Thus, $k_0\in P^-$. Since $P^-\cap K\subset Z$, the claim follows.

\vs\noindent\underline{Step 2}: We claim that $d(Ka_ik_0a_i^{-1},K)\to 0$ as $i\to\infty$.

\vs Write $k_i=k_0l_i$ where $l_i\in K$ and $l_i\to e$ as $i\to\infty$. Then
\begin{equation}\label{eq_st0}
d(Ka_ik_i,Ka_i)=d(Ka_ik_0a_i^{-1},Ka_il_i^{-1}a_i^{-1})\to 0\;\;\hbox{as}\;\; i\to\infty.
\end{equation}
Thus, the sequence $\{a_il_i^{-1}a_i^{-1}\}$ is bounded, and hence, we may assume that it converges.
Since $l_i^{-1}\to e$ as $i\to\infty$ and $\mathfrak{g}=\mathfrak{u}^-\oplus\mathfrak{z}\oplus\mathfrak{u}^+$, we obtain that
$l_i^{-1}=u^-_iz_iu^+_i$ for some $u_i^-\in U^-$, $z_i\in Z$, and $u_i^+\in U^+$ such
that $u_i^-\to e$, $z_i\to e$, and $u_i^+\to e$ as $i\to\infty$. It follows from (\ref{eq_uzu}) that
$a_iu^-_ia_i^{-1}\to e$ and $a_iz_ia_i^{-1}\to e$ as $i\to\infty$. Hence, $a_il_i^{-1}a_i^{-1}\to u^+$ as $i\to\infty$ for some $u^+\in U^+$.
On the other hand, by Step 1, passing to a subsequence,
we get $a_ik_0a_i^{-1}\to z$ as $i\to\infty$ for some $z\in Z$. Thus, by (\ref{eq_st0}), $z^{-1}u^+\in K$ as $i\to\infty$. Since $P^+\cap K\subset Z$,
we deduce that $u^+=e$. Hence, $a_il_i^{-1}a_i \to e$ as $i\to\infty$, and the claim follows from (\ref{eq_st0}).

\vs\noindent\underline{Step 3}: We claim that $a_i k_0a_i^{-1}=d_ik_0d_i^{-1}$ for some bounded sequence $\{d_i\}\subset \mathcal{C}$.

\vs Recall that the system of simple roots $\Pi$ is a basis of the dual space of the Lie algebra of $A$.
Hence, we may write $a_i=b_ic_i$ where $b_i,c_i\in A^+$ such that $\alpha(\log b_i)=0$
for every $\alpha\in \Pi\cap\Phi_+$ and $\alpha(\log c_i)=0$ for every $\alpha\in\Pi-\Phi_+$.
Then $c_i$ commutes with $Z$, and since $\alpha(b_i)=\alpha(a_i)$ for every $\alpha\in\Pi-\Phi_+$,
the sequence $\{b_i\}$ is bounded.
Let $c_0\in A^+$ be such that $\alpha(\log c_0)=C$ for every $\alpha\in \Pi\cap \Phi_+$ and $\alpha(\log c_0)=0$ for every $\alpha\in\Pi-\Phi_+$.
Then $d_i=b_ic_0\in \mathcal{C}$ and $a_ik_0a_i^{-1}=d_ik_0d_i^{-1}$ as required.

\vs Taking a subsequence, we obtain that $d_i\to
d_0$ as $i\to\infty$ for some $d_0\in \mathcal{C}$ and $d_0k_0d_0^{-1}\in K$.
This implies that $d_0k_0=kd_0$ for some $k\in K$. Since $d_0$ lies in the interior of $A^+$, it
follows from the uniqueness properties of the Cartan decomposition
(Lemma \ref{unique}) that $k_0=k\in M$. This is a contradiction.
\end{proof}

\begin{Thm}[\bf The strong wavefront lemma] \label{l_wave}
 Let $\Cal C$ be a closed subset contained in $A^+$ with a positive distance
 from the walls of $A^+$.
Then for any neighborhoods $U$ of $e$ in $K$ and $V$ of $e$ in $A$, there
 exists a neighborhood $\Cal O$ of $e$ in $G$ such that
 $$\Cal O g \cup g\Cal O \subset (k_1U)(a VM)(k_2U) $$
 for any $g=k_1ak_2\in K\Cal C K$. \end{Thm}

\begin{proof}
Without loss of generality, we may assume that $\Cal C=A^+(C)$ for
some $C>0$, where
$$
A^+(C)=\{a\in A^+: \alpha(\log a)\ge C\;\hbox{ for all $\alpha\in\Phi^+$}\}.
$$
Replacing $V$  by a smaller subset if necessary, we
may assume that $A^+(C)V$ is contained in the interior of
$A^+(C/2)$.

\vs\vs\noindent\underline{Step 1}: We claim that there exists an $\e$-neighborhood
$\Cal O$ of $e$ in $G$ such that
$$\Cal Og \subset Ka VK \quad\text{for all $g=k_1ak_2\in KA^+(C)K$}.$$
We take $\varepsilon>0$ to be sufficiently small such that the
$\e$-neighborhood of $e$ in $A$ is contained in $V$. Let
$\mathcal{O}$ be the $\e$-neighborhood of $e$ in $G$. 
For $h\in \mathcal{O}a$, write $h=k_1' b
k_2'\in KA^+K$. Then $d(Kbk_2'a^{-1},K)=d(Kbk_2',Ka)<\varepsilon$.
Since $d(Kb,Ka)\le d(Kbk_2',Ka)$ by Lemma \ref{l_hel}, we have
\begin{equation}\label{s1}
d(Kb,Ka)<\varepsilon\end{equation}
 and hence, $h=k_1'bk_2'\in
KVaK$. This shows that $\Cal O a\subset KaVK$. Since $\Cal O$ is
$K$-invariant,
 $$\Cal O k_1 a k_2\subset \Cal O a k_2\subset KaV K$$
This proves the claim.

\vs\vs\noindent\underline{Step 2}: We claim that there exists an $\e$-neighborhood $\Cal O$ of $e$
in $G$ such that
\begin{equation}\label{s2}\Cal Og \subset
KA^+ (k_2U) \quad\text{for all $g=k_1ak_2\in
KA^+(C)K$}.\end{equation}

 Let $U_0$ be a neighborhood of $e$ in $K$ such that
$kU_0k^{-1}\subset U $ for all $k\in K$. Choose $\e>0$ that satisfies Step 1 above,
and Proposition \ref{strong} holds with respect to $2\e$, $U_0$, and $A^+(C/2)$.
Take $\Cal O$ to be the $\e$-neighborhood of $e$ in
$G$. Then for any $h=k_1'bk_2'\in \Cal O a$, we have by
(\ref{s1}),
$$d(Kbk_2', Kb)\le d(Kb k_2', Ka)+ d(Kb,Ka)<2\varepsilon.
$$
Since $\Cal O a\subset KA^+(C/2)K$, by Proposition \ref{strong},
this implies that $k_2'\in MU_0$, and hence $\Cal O a\subset
KA^+MU_0$. Thus,
$$\Cal O (k_1 ak_2)\subset KA^+MU_0k_2\subset KA^+M(k_2U) .$$

\vs\vs\noindent\underline{Step 3}: We claim that there exists a neighborhood $\Cal O$ of
$e$ in $G$ such that
$$\Cal O g \subset (k_1 U) A^+K \quad\text{for all $g=k_1ak_2\in KA^+(C)K$}.$$

Recall that by the Iwasawa decomposition, $G=KAN$, where $N$ is the
subgroup of $G$ corresponding to the sum of the positive root
spaces. Since the Weyl group acts transitively on the sets
of simple roots, there exists an element $w\in K$ that normalizes $A$ such that
$wA^+(C) w^{-1}={(A^+ (C))}^{-1}$.

Let $U_0$ be a neighborhood of $e$ in $K$ such that $U_0(wU_0^{-1}
w^{-1}) \subset U$.
Let $\Cal O$ be the neighborhood from Step $2$ with respect to $U_0$ and ${\Cal O}_1$
a neighborhood of $e$ in $AN$ such that $ c^{-1} {\Cal O_1} c
\subset  {\Cal O}$ for all $c\in A^+$.
Since $w\in K$, $w{\Cal O}w^{-1}={\Cal O}$.
Conjugating (\ref{s2}) (with respect to $U_0$)
by $w$, we have for any $b\in w A^+(C) w^{-1}$,
\begin{equation*}
\mathcal{O} b\subset K {(A^+)}^{-1}(w U_0 w^{-1}).
\end{equation*}
For $a\in A^+(C)$, $a^{-1}\in w A^+(C)w^{-1}$ and hence
\begin{equation}\label{eq_wave3}
a^{-1} {\mathcal{O}_1}\subset {\mathcal{O}} a^{-1} \subset K
{(A^+)}^{-1} (w U_0w^{-1}).
\end{equation}
By taking the inverse of (\ref{eq_wave3}),
\begin{equation}\label{ss3} {\Cal O}_1 ^{-1} a\subset (wU_0^{-1} w^{-1}) A^+
K .\end{equation}

Since the product map $K\times AN \to G$ is a diffeomorphism,
$U_0\Cal O_1 ^{-1}$ is a neighborhood of $e$ in $G$. Therefore,
there exists a neighborhood $\Cal O_2$ of $e$ in $G$ such that
$k^{-1}\Cal O_2 k\subset U_0 \Cal O_1^{-1}$ for all $k\in K$. Then
by (\ref{ss3}), $${\Cal O}_2k_1 ak_2\subset k_1 U_0 \Cal O_1^{-1}
a k_2\subset k_1 U_0
 (wU_0^{-1} w^{-1}) A^+ K
\subset k_1 UA^+K .$$ This finishes the proof of Step 3.

\vs\vs

By the above three steps, we obtain a neighborhood $\Cal O_1$ of
$e$ in $G$ such that for all $g=k_1ak_2\in KA^+(C)K$,
$$\Cal O_1 g\subset
K aVK\cap KA^+(k_2 U)\cap (k_1U)A^+K .$$  By the
uniqueness the Cartan decomposition (Lemma \ref{unique}),
$$K aVK\cap KA^+(k_2 U)\cap (k_1U)A^+K =(k_1U)(aVM)(k_2U) .$$
Hence, we have shown that for every neighborhood $U$ of $e$ in $K$
and every neighborhood $V$ of $e$ in $A$, there exists a neighborhood $\mathcal{O}_1$ of $e$ in $G$
such that
\begin{equation}\label{w1}
\Cal O_1 g\subset (k_1U)(aVM)(k_2U)\quad\hbox{for all $g=k_1ak_2\in KA^+(C)K$}.
\end{equation}

Let $U_0$ be a neighborhood of $e$ in $K$ such that $kU_0^{-1}k^{-1}\subset U$ for every $k\in K$
and $V_0$ a neighborhood of $e$ in $A$ such that $wV_0^{-1}w^{-1}\subset V$, where $w$ is the element of $K$ defined in Step 3.
Let $\Cal O_2$ be a neighborhood of $e$ in $G$ such that 
$$
\Cal O_2 g\subset (k_1U_0) (aV_0M)(k_2U_0)\quad\hbox{for all $g=k_1ak_2\in KA^+(C)K$}.
$$
By taking the inverse, we have 
$$
g\Cal O_2^{-1}\subset (U_0^{-1}k_1) (aV_0^{-1}M)(U_0^{-1}k_2)\quad\hbox{for all $g=k_1ak_2\in KA^+(C)^{-1}K$}.
$$
Conjugating by $w$, we get
$$
g(w\Cal O_2^{-1}w^{-1})\subset (k_1U) (aVM)(k_2U)\quad\hbox{for all $g=k_1ak_2\in KA^+(C)K$}.
$$
Set $\Cal O=\Cal O_1\cap (w\Cal O_2^{-1}w^{-1})$ to finish the proof.
\end{proof}

\section{Contractions on $G/B$}
\label{boundary}

Let $G$ be a connected noncompact semisimple Lie group with finite center, $K$ a maximal compact subgroup of $G$,
and $B$ a maximal connected split solvable subgroup of $G$.
Let $N$ be the unipotent radical of $B$. The normalizer $P$ of $N$ in $G$  is the unique minimal parabolic subgroup of $G$ containing $B$.

\begin{Lem}\label{langlands}
There exists a split Cartan subgroup $A$ of
$G$ with respect to $K$ and an ordering on the root system
$\Phi$ such that
\begin{itemize}
\item $G=KA^+K$.
\item $B=AN$.
\item $N$ is the subgroup generated by all positive root subgroups of $G$ with respect to $A$.
\item $P=MAN$, where $M$ is the centralizer of $A$ in $K$ and $M=K\cap P$.
\end{itemize}
\end{Lem}
\begin{proof} Take any split Cartan subgroup $A_0$ and a Weyl chamber
$A_0^+$ such that the Cartan
decomposition $G=KA_0^+K$ holds. Set $P_0=M_0A_0N_0$, where $M_0$ is the centralizer of $A_0$ in $K$, and
$N_0$ is
the subgroup generated by all positive root subgroups of $G$ with respect to $A_0$.
Note that $P_0$ is a minimal parabolic subgroup of $G$ (see \cite[Sec.~1.2.3]{wa}).
By the Iwasawa decomposition, $G=KP_0$. Since all minimal parabolic subgroups are conjugate to each other ([loc.cit.]),
there exists $k\in K$
 such that $P=kP_0k^{-1}$. Set $A=kA_0k^{-1}$, $A^+=kA^+_0k^{-1}$, and
$M=kM_0k^{-1}$. It is clear that $N=kN_0k^{-1}$. Since $AN$ is normal in $P$, it is the unique maximal connected
split solvable subgroup in $P$. Thus, $B=AN$.
\end{proof}

Let $A$ be a split Cartan subgroup as in Lemma \ref{langlands}.
Denote by $M$ be the centralizer of $A$ in $K$ and by $N^-$
the subgroup generated by all negative root subgroups of $G$. One can check that the map
$$
N^-\times M\to G/B: (n,m)\mapsto nmB
$$
is a diffeomorphism onto its image. Moreover, it follows from the properties of Bruhat decomposition
that $G/B-\pi(N^-\times M)$ is a finite union of closed submanifolds of smaller dimensions. On the other hand, by
the Iwasawa decomposition, the map
$$
K\to G/B: k\mapsto kB.
$$
is a diffeomorphism. Thus, we have a map $N^-\times M\to K$. Since $M$ normalizes $B$, this map is right $M$-equivariant.
Denote by $S\subset K$ the submanifold which is the image of $N^-\times \{e\}$ under this map.

Let $\nu$ be the probability Haar measure on $K$ and $\tau$ the probability Haar measure on $M$.
Then $\nu=\sigma\otimes \tau$ for some finite smooth measure $\sigma$ on $S$ (see \cite[p.~73]{wa}),
i.e.,
\begin{equation}\label{eq_sm}
\int_K f\,d\nu =\int_M\int_S f(sm)\,d\sigma(s)d\tau(m),\quad f\in C(K).
\end{equation}

For $C>0$, put
$$A^+(C)=\{a\in A: \alpha(\log a) \ge  C\;\;\hbox{ for all $\alpha\in\Phi^+$} \} .$$

\begin{Lem}\label{k1}
Let $\Omega$ be a Borel subset of $K$, $s\in S$, $m\in M$, and $k\in K$. For $a\in A$, $r\in K$, and $U\subset K$, define
\begin{equation*}
\Omega_{r}(k,a,U)=\{l\in \Omega:  lkB\in a^{-1}r^{-1}UB\}.
\end{equation*}
\begin{enumerate}
\item Let $VW\subset K$ be an open neighborhood of $s$ where $V\subset S$ and $W\subset M$ are open subsets.
Then for every $\varepsilon>0$, there exists
$C>0$ such that $$\nu((\Omega\cap m^{-1}SWk^{-1})-\Omega_{sm}(k,a,VW))<\varepsilon\;\;\hbox{ for all $a\in A^+(C)$.}$$
\item Let $U\subset K$ be a Borel subset such that $s\notin\overline{UM}$. Then for every $\varepsilon>0$, there exists
$C>0$ such that $$\nu(\Omega_{sm}(k,a,U))<\varepsilon\;\;\hbox{for all $a\in A^+(C)$.}$$
\end{enumerate}
\end{Lem}

\begin{proof}
Note that $\Omega_{sm}(k,a,U)=m^{-1}\cdot(m\Omega)_{s}(k,a,U)$. Hence, replacing $\Omega$ by $m\Omega$, we may assume that $m=e$.
Also, $\Omega_{s}(k,a,U)=(\Omega k)_s(e,a,U)\cdot k^{-1}$. Thus, we may assume that $k=e$.

The proof is based on the observation that elements in the interior of $A^+$ act on $N^-B$
as contractions (cf. \cite[Proposition~8.2.5]{Zi}).
Denote by $\pi$ the map $\pi:N^{-}\to G/B:n\mapsto nB$.
Note that for $a\in A$ and $V_0\subset S$, 
\begin{equation}\label{eq_D}
\pi^{-1}(a^{-1}V_0B)=a^{-1}\pi^{-1}(V_0B)a.
\end{equation}
Let $D$ be a compact set in $S$ such that $\sigma(S-D)<\e/2$.

Suppose that $s\in VW$. Let $W_0\subset W$ be a compact set such that $\tau({W-W_0})<\e/2$.
Then $s^{-1}VW$ is a neighborhood of $W_0$. By uniform continuity, there exists a neighborhood $V_0$ of $e$ in $S$
such that $V_0W_0\subset s^{-1}VW$. For each $a\in A^+$, the map $N^-\to N^-: x\to a^{-1}x a$ expands any neighborhood
of $e$ at least by the
factor of $\min_{\alpha\in \Phi^+} e^{\alpha(\log a)}$.
Hence, there exists $C>0$ such that $\pi^{-1}(DB) \subset a^{-1} \pi^{-1}(V_0B) a$ for all $a\in A^+(C).$
Then by (\ref{eq_D}), $DB\subset a^{-1}V_0B$. It follows that $DW_0B\subset a^{-1}V_0W_0B$ and
\begin{align*}
\nu(\Omega_{s}(e,a,VW))&\ge \nu(\Omega\cap DW_0)\ge \nu(\Omega\cap SW)-\nu((S-D)M)-\nu(S(W-W_0))\\
&\ge \nu(\Omega\cap SW)-\e.
\end{align*}
This proves the first part of the lemma. 

To prove the second part, we observe that there exists an open subset $V\subset S$ such that $s\in V\subset K-\overline{UM}$.
Since $\nu(m^{-1}SMk^{-1})=1$, it suffices to apply the first part with $W=M$.
\end{proof}

\begin{Lem}\label{l_v}
Let $W\subset M$ and $\Omega\subset K$ be Borel subsets. Then every $k\in K$,
$$
\int_M \nu(\Omega\cap m^{-1}SW k^{-1})\,d\tau(m)=\nu(\Omega)\tau(W).
$$
\end{Lem}

\begin{proof}
Without loss of generality, $k=e$.
Since $M$ normalizes $N^-$, we have $m^{-1}Sm\subset S$. Since $\nu(K-SM)=0$, we may assume that $\Omega\subset SM$. 
Moreover, it suffices to prove the lemma for $\Omega= \Omega_S\Omega_M$ with Borel sets $\Omega_S\subset S$ and $\Omega_M\subset M$.
By (\ref{eq_sm}), 
$$
\nu(\Omega\cap m^{-1}SW)=\sigma(\Omega_S)\tau(\Omega_M\cap m^{-1}W).
$$
For $Y\subset M$, denote by $\chi_Y$ the characteristic function of $Y$. We have
\begin{align*}
&\int_M \tau(\Omega_M\cap m^{-1}W)\,d\tau(m)=\int_M \int_M \chi_{\Omega_M}(l) \chi_{m^{-1}W}(l)\, d\tau(l)d\tau(m)\\
&=\int_M \chi_{\Omega_M}(l) \left(\int_M \chi_{Wl^{-1}}(m)\,d\tau(m)\right) d\tau(l)=\tau(\Omega_M)\tau(W).
\end{align*}
Since $\nu(\Omega)=\sigma(\Omega_S)\tau(\Omega_M)$, this proves the lemma.
\end{proof}

\section{Volume estimates}\label{s_vol}

Let $G$ be a connected noncompact semisimple Lie group with finite center and $G=KA^+K$ a Cartan decomposition.
Let $dk$ denote the probability Haar measure on $K$, $dt$ the
Lebesgue measure on the Lie algebra $\mathfrak{a}$ of $A$, and $da$ the
Haar measure on $A$ derived from $dt$ via the exponential map.
We denote by $m$ the Haar measure
on $G$ which is normalized so that for any $f\in C_c(G)$,
\begin{equation}\label{cartan}\int_G f\, dm
=\int_{K}\int_{A^+}\int_K f(k_1a k_2) \xi(\log a) dk_1 da
dk_2,\end{equation}
where
\begin{equation}\label{eq_xi}
\xi(t)=\prod_{\alpha\in \Phi^+}\sinh (\alpha(t)) ^{m_\alpha},\;\; t\in\mathfrak{a},
\end{equation}
and $m_\alpha$ is the dimension of the root subspace corresponding to $\alpha$.
In particular, for any measurable subset $D\subset A^+$, we
have 
\begin{equation}\label{eq_B}
\vol (KDK)=\int_{D}\xi(\log a)\, da.
\end{equation}

Let $\|\cdot \|$ be a Euclidean norm on $\frak a$,
that is, for some basis $v_1, \ldots, v_n$ of $\frak a$, $\|\sum_i c_i v_i\|
=\sqrt {\sum_{i} c_i^2}$.
We assume that $\|\cdot\|$ is invariant under the Weyl group action.
For instance, $\|\cdot \|$ can be taken to be the norm
induced from an invariant Riemannian metric $d$ on the
symmetric space $K\ba G$, i.e.,
 $$\|t\|=d (K,K\exp(t))\;\;\hbox{ for } t\in \mathfrak{a}.$$
For $T>0$ and $\mathfrak{t}\subset \mathfrak{a}$, set
$$
\mathfrak{t}_T=\{t\in \mathfrak{t}:\; \|t\|< T  \}.
$$
Let $\rho=\frac{1}{2}\sum_{\alpha\in\Phi^+}m_\alpha\alpha$. 
One can check that the maximum of $\rho$ on $\bar{\mathfrak{a}}_1$ is
achieved at a unique point contained in the interior of $\mathfrak{a}^+$, which we call the {\it barycenter} of $\mathfrak{a}^+$.

We present a simple derivation of the asymptotics of the volume of Riemannian balls in a symmetric space of noncompact type
(see also \cite[Theorem 6.2]{kni}).

\begin{Lem}\label{p_vol_as0}
Let $\mathfrak{q}$ be a convex cone in $\mathfrak{a}^+$ centered at the origin that contains the barycenter of $\mathfrak{a}^+$
in its interior. Then for some $C>0$ (independent of $\mathfrak{q}$),
$$
\int_{\mathfrak{q}_T} \xi(t)\, dt \sim_{T\to\infty} C\cdot T^{(r-1)/2} e^{\delta T},
$$
where $r=\operatorname{\mathbb{R}-rank}(G)=\dim A$ and $\delta=\max\{2\rho(t):t\in \bar{\mathfrak{a}}_1\}$. In particular,
$$
\int_{\mathfrak{q}_T} \xi(t)\, dt \sim_{T\to\infty} \int_{\mathfrak{a}^+_T} \xi(t)\, dt.
$$
\end{Lem}
\begin{proof}
We have
\begin{align*}
\int_{{\mathfrak{q}_T}} \xi (t)\, dt
=\frac{1}{2^{|\Phi^+|}}\int_{{\mathfrak{q}_T}} e^{2\rho(t)}\, dt + \hbox{``other
terms'',}
\end{align*}
where the ``other terms'' are linear combinations of integrals of the
form
$
\int_{{\mathfrak{q}_T}} e^{\lambda(t)}\, dt
$
such that $2\rho-\lambda=\sum_{\alpha\in\Phi^+} n_\alpha \alpha$ for some $n_\alpha\ge 0$. 
In particular, $2\rho>\lambda$ in the interior of the
Weyl chamber $\mathfrak{a}^+$.
Since the maximum of
$2\rho$ in ${\bar{\mathfrak{a}}_1}$ is achieved in the interior of $\mathfrak{q}$,
$$
\max\{2\rho(t): t\in {\bar{\mathfrak{q}}_1}\}>\max\{\lambda (t): t\in
{\bar{\mathfrak{q}}_1}\}\stackrel{def}{=}\delta'.
$$
Thus,
$$
\int_{{\mathfrak{q}_T}} e^{\lambda(t)}\, dt\le e^{\delta'
T}\vol({\mathfrak{q}_T})\ll e^{\delta'
T}T^r =o(e^{\delta T})
$$
as $T\to\infty$.
It remains to show that for some $C>0$ independent of $\mathfrak{q}$,
\begin{equation}\label{eq_dom}
\int_{{\mathfrak{q}_T}} e^{2\rho(t)}\, dt\sim_{T\to\infty} C\cdot T^{(r-1)/2}
e^{\delta T}.
\end{equation}
Making a change variables and decomposing ${\frak{q}_1}$ into slices parallel to the
hyperplane $\{2\rho=0\}$, we have
$$
\int_{{{\frak{q}}_T}} e^{2\rho(t)}\, dt = T^r
\int_{{{\frak{q}}_1}} e^{2T\rho(t)}\, dt= T^r
\int_0^\delta
e^{Tx}\phi(x)\, dx,
$$
where $\phi(x)=\vol_{r-1}({\mathfrak{q}}_1\cap \{2\rho=x\})$.

First, we show that for some $c_1>0$ independent of $\mathfrak{q}$,
\begin{equation}\label{eq_phi_a} \phi(x)\sim_{x\to\delta^-} c_1\cdot (\delta -
x)^{(r-1)/2}.
\end{equation}
We identify $\mathfrak a$ with
the set $\{(t_1, \cdots, t_{r})\in \br^r\}$
and denote by $Q$ the positive quadratic form on $\mathfrak a$
defined by the norm.
After a linear change of variables, we may assume that $2\rho(t)=t_r$ and
$$
Q(t_1,\ldots,t_r)=\sum_{i=1}^{r-1} \alpha_i (t_i-\beta_i
t_r)^2 +\alpha_r t_r^2
$$
for some $\alpha_i>0$ and $\beta_i\in\mathbb{R}$.
It is clear that
the maximum of $2\rho(t)$ on the set
$\bar{\mathfrak{a}}_1=\{t:Q(t)\le 1\}$ is achieved
when $t_i=\beta_i t_r$ for $i=1,\ldots, r-1$ and $\alpha_r t_r^2=1$.
This implies that $\delta=\alpha_r^{-1/2}$. Since for $x$ close to
$\delta$, the set ${\frak{q}}^+_1\cap \{2\rho=x\}$ is defined by
the condition $Q(t_1,\ldots,t_{r-1},x)\le 1$, we have
\begin{align*}
\phi(x)&=\vol_{r-1}
\left(\left\{(t_1,\ldots,t_{r-1}):\sum_{i=1}^{r-1} \alpha_i (t_i-\beta_i x)^2 \le 1-
\alpha_r x^2 \right\}\right)\\
&= c_1\cdot (1-\alpha_r x^2)^{(r-1)/2}=c_1\cdot
(1-x^2/\delta^2)^{(r-1)/2}
\end{align*}
for a constant $c_1>0$.
This proves (\ref{eq_phi_a}). Now by (\ref{eq_phi_a}) and the L'Hopital's rule,
$$\beta(x)\stackrel{def}{=}\int_0^x \phi(\delta-u)\,du \sim_{x\to
0^+} c_2\cdot x^{(r+1)/2}$$
for some $c_2>0$. Thus, by the Abelian theorem (see \cite[Corollary~1a, p.~182]{wid}),
\begin{align*}
\int_0^\delta e^{Tx}\phi(x)\, dx &= e^{\delta T} \int_0^\delta e^{-Tx}\phi(\delta-x)dx
=  e^{\delta T} \int_0^\infty e^{-Tx}d\beta(x)\\
&
\sim_{T\to\infty} c_3\cdot T^{-(r+1)/2} e^{\delta T}
\end{align*}
for some $c_3>0$. It is clear that $c_3$ is independent of $\frak q$.
This proves (\ref{eq_dom}) and the proposition.
\end{proof}

For $T>0$ and $R\subset G$, define
\begin{equation}\label{eq_rt}
R_T=\{r\in R:\; d(K,Kr)< T  \}.
\end{equation}
Since $G_T=K{A_T^+}K$, combining the previous lemma and (\ref{eq_B}), we deduce

\begin{Cor}\label{p_vol_as}
For some $C>0$,
$$
\vol (G_T)\sim_{T\to\infty} C\cdot T^{(r-1)/2} e^{\delta T}.
$$
\end{Cor}
In particular, we have
\begin{Lem} \label{vol} For some
functions $a(\e)$ and $b(\e)$ such that $a(\e)\to 1 $ and $b(\e)\to 1$ as $\e\to 0^+$,
we have
$$  a(\e) \le \liminf_{T\to\infty}\frac{\vol (G_{T-\e})}{\vol(G_{T})}
  \le  \limsup_{T\to\infty}\frac{\vol(G_{T+\e})}{\vol (G_T)}\le b(\e) .$$
\end{Lem}

For $T,C>0$ and $R\subset A$, define
\begin{eqnarray}\label{eq_aaa}
R(C)&=&\{r\in R:\; \alpha(\log r) \ge C \;\; \text{ for
all $\alpha\in \Phi^+ $} \},\\
R_T(C)&=&R_T\cap R(C).\nonumber
\end{eqnarray}

\begin{Lem}\label{cone}
Let $Q\subset A^+$ be a convex cone centered at the origin that contains the barycenter in its interior.
Then for any fixed $C>0$,
$$
\int_{Q_T(C)} \xi(\log a)\, da\sim_{T\to\infty}\int_{A^+_T} \xi(\log a)\, da.
$$
In particular,
$$\vol (KA_T^+K)\sim_{T\to\infty} \vol (KA_T^+(C) K) .$$
\end{Lem}
\begin{proof}
It suffices to prove the lemma when $Q$ is contained in the interior of $A^+$.
Then there exists $T_0>0$ such that $Q_T(C)\supset Q_T-Q_{T_0}$ for all sufficiently large $T>0$. Thus, the lemma follows from Lemma \ref{p_vol_as0}.
\end{proof}

\section{Equidistribution of solvable flows}\label{shahsection}

Let $G$ be a connected noncompact semisimple Lie group with
finite center which is realized as a closed subgroup of a Lie group $L$,
$\Lambda$ a lattice in $L$, and $Q$ a closed subgroup of $G$
that contains a maximal connected split solvable subgroup of $G$.
In this section we investigate the distribution of orbits of $Q$ in the homogeneous space $\Lambda\ba L$.

Let $K_0$ be a maximal compact subgroup of $G$ and $A_0$ a
split Cartan subgroup of $G$ with respect to $K_0$. Then $G=K_0A_0^+K_0$
for
any positive Weyl chamber $A_0^+$ in $A_0$.
Denote by $d$ an invariant Riemannian metric on $K_0\ba G$.
For $g\in G$, $R\subset G$, $\Omega\subset K_0$, and $T>0$, define
$$R_T(g, \Omega)= \{r\in R: d(K_0, K_0g r) <T, \,r\in  g^{-1}K_0 {A_0^+} \Omega \}$$
and $R_T(g)=R_T(g, K_0)$.
Note that $$R_T(k_0g,\Omega)=R_T(g,\Omega)=g^{-1}(gR)_T(e,\Omega)$$ for any $k_0\in K_0$.

The main result in this section is Theorem \ref{Peq} on equidistribution of the
sets $Q_T(g, \Omega)$ as $T\to \infty$.
Let $\mu$ be the Haar measure on $L$ such that $\mu(\Lambda\ba L)=1$,
$\nu_0$ the probability Haar measure on $K_0$, and $\rho$ a right invariant Haar measure on $Q$.
Denote by $G_n$ the product of all noncompact simple factors of $G$.

\begin{Thm}\label{Peq}
Suppose that for $y\in \Lambda \ba L$, the orbit $yG_n$ is dense in $\Lambda\ba L$.
Then for $g\in G$, any Borel subset $\Omega\subset K_0$
with boundary of measure zero
and $f\in C_c(\Lambda\ba L)$,
\begin{equation*}\lim_{T\to \infty} \frac{1}{\rho(Q_T(g))}\int_{Q_T(g,
\Omega)} f(y q^{-1})\, d\rho(q)= \nu_0(M_0\Omega)
\int_{\Lambda \ba L}f \, d\mu.
\end{equation*}
\end{Thm}

The rest of this section  is devoted to a proof of Theorem
\ref{Peq}. We start by stating a theorem of Shah (\cite[Corollary~1.2]{Sh}).
Recall that a sequence $\{g_i\}\subset G$ is called {\it strongly divergent} if 
for every projection $\pi$ from $G$ to its noncompact factor, $\pi(g_i)\to\infty$ as $i\to\infty$.

\begin{Thm}[\bf Shah]\label{shah}
Suppose that for $y\in \Lambda \ba L$, the orbit $yG_n$ is dense in $\Lambda\ba L$.
Let $\{g_i\}\subset G$ be a strongly divergent sequence.
Then for any $f\in C_c(\Lambda\ba L)$ and any Borel subset $U$ in $K_0$ with boundary of measure zero, 
\begin{equation*} \lim_{i\to\infty} \int_{U}
f(y kg_i )\, d\nu_0(k) =\nu_0 (U) \int_{\Lambda \ba L}f \, d\mu.
\end{equation*}
\end{Thm}

\begin{Rmk} {\rm
Although this theorem was stated in \cite{Sh} only for the case $U=K_0$, its proof works equally well
when $U$ is an open subset of $K_0$ with boundary of measure zero, and approximating Borel sets by open sets, one can check
that Theorem \ref{shah} holds in the above generality.
}
\end{Rmk}

Let $K=g^{-1} K_0 g$ and $B$ be a maximal connected split solvable subgroup of $G$.
By Lemma \ref{langlands}, there exists a a split Cartan subgroup $A$ with
an ordering on the root system $\Phi$ such that $G=KA^+ K$ and $B=AN$, where
$N$ is the subgroup generated by all positive root subgroups of $A$ in $G$.
Let $M$ be the centralizer of $A$ in $K$. 

Let $m$ be the Haar
measure on $G$ such that (\ref{cartan}) holds with respect to
$K_0$ and $A_0^+$.
It follows from the uniqueness of the Haar measure that for every $g\in G$, there exists $c_g>0$ such that
\begin{equation}\label{kp}
\int_G f\, dm = c_g\int_{K_0}\int_{B} f(g^{-1}kg b)\,d\rho(b)d\nu_0(k),\;\; f\in C_c(G).
\end{equation}
In particular,
\begin{equation}\label{eq_rho_vol}
c_g\rho(B_T(g))=\vol (G_T(g))=\vol (G_T(e)).
\end{equation}
We normalize $\rho$ so that $c_e=1$.
Let $\nu$ be the probability Haar measure on $K$.

We use notations from Section \ref{boundary}. In particular, (\ref{eq_sm}) holds.

\begin{Prop}\label{weakc}
Suppose that for $y\in \Lambda \ba L$, the orbit $yG_n$ is dense in $\Lambda\ba L$.
Let $U=VW$ be an open neighborhood of $e$ in $K$ where $V$ is an open neighborhood of $e$ in $S$, and $W$ is an open neighborhood of $e$ in $M$ such that
$\sigma(\partial V)=\tau(\partial W)=0$
and $\Omega$ a Borel subset of $K_0$
such that $\nu_0(\partial\Omega)=0$.
Then for any $f\in C_c(\Lambda\ba L)$,
 $$\frac{1}{\rho(B_T(g))} \int_{U}\int_{B_T(g,\Omega)}
f(y b^{-1} k^{-1})\, d\rho(b) d k \to \nu({U})
\nu_0({M_0\Omega}) \int_{\Lambda\ba L} f \, d\mu$$ as $T\to
\infty$.
\end{Prop}
\begin{proof}
Since both $A_0$ and $gAg^{-1}$ are split Cartan subgroups
with respect to $K_0$, there exists $k_0\in K_0$ such that $A_0=k_0gAg^{-1}k_0^{-1}$.
Hence, replacing $g$ by $k_0g$, we may assume that $A=g^{-1}A_0g$.

Let $d$ be the Riemannian metric on $K\ba G$ induced from the metric on $K_0\ba G$ by
the map $Kx\mapsto Kgx$, and notations $A_T^+$, $A_T^+(C)$, and $A^+(C)$ be
defined as in Section \ref{s_vol}.
For instance, $$A_T^+=\{a\in A^+: d(K_0g, K_0ga)< T\} .$$
Since $KA^+_TK=g^{-1}G_T(e)g$, it follows from (\ref{eq_B}) and (\ref{eq_rho_vol}) that
for every $T>0$,
\begin{equation}\label{eq_ar0}
\rho(B_T(e))=\int_{A_T^+} \xi(\log a)\,da
\end{equation}
and by Lemma \ref{cone},
\begin{equation}\label{eq_ar}
\int_{A_T^+-A_T^+(C)} \xi(\log a)\,da= o(\rho(B_T(e)))\quad\hbox{as $T\to\infty$.}
\end{equation}

Since $B_T(g, \Omega)=B_T(g, M_0\Omega)$, we may assume
without loss of generality that $M_0\Omega =\Omega$.

We have
\begin{eqnarray*}
UB_T(g, \Omega)&=&g^{-1}K_0{A_{0\,T}^+}(e)\Omega\cap UB=KA_T^+(g^{-1}\Omega g)g^{-1}\cap UB\\
&=&\{k_1 ak_2g^{-1}: k_1\in K,
a\in A^+_T, k_2\in \Omega_{k_1}(a)\},
\end{eqnarray*}
where
\begin{eqnarray*}
\Omega_{{k}_1}(a)=\{{k}_2\in
g^{-1}\Omega g:  {k}_2g^{-1}B\in {a}^{-1} {k}_1^{-1}UB\}.
\end{eqnarray*}
By (\ref{cartan}) and (\ref{kp}),
\begin{align}\label{sum}
&c_g\int_{U}\int_{B_T(g, \Omega)} f(y b^{-1} k^{-1})\, d\rho(b) d\nu(k)=\int_{UB_T(g,\Omega)}f(y x^{-1})\, dm(x)\\
 &=\nonumber \int_{{k}_1{a}{k}_2g^{-1}\in U
B_T(g, \Omega) } f(yg{k}_2^{-1}{a}^{-1} {k}_1^{-1})\,
\xi(\log {a}) d\nu({k}_2) d{a} d\nu({k}_1)\\ &=\nonumber \int_{{k}_1\in {K}}
\int_{{a}\in
 A_T^+}\int_{{k}_2 \in
\Omega_{{k}_1}({a}) } f(yg{k}_2^{-1}{a}^{-1} {k}_1^{-1})\,
\xi(\log{a}) d\nu({k}_2) d{a} d\nu({k}_1).
\end{align}

\vs\vs\noindent\underline{Step 1}: We claim that
for every $m\in M$,
\begin{multline}\label{ddd1}
\lim_{T\to \infty} \frac{1}{\rho(B_T(e))} \int_{V}
\int_{A^+_T} \int_{\Omega_{sm}({a}) }
f(yg {{k}}^{-1}{a}^{-1} (sm)^{-1})\, \xi(\log {a}) d\nu({k}) d{a}
d\sigma(s)\\ =\sigma(V)\cdot\nu({g^{-1}\Omega g\cap m^{-1}SWk_g^{-1}})\cdot \int_{\Lambda\backslash L}f\, d\mu,
 \end{multline}
where
$k_g\in K$ such that $k_gB=g^{-1}B$.
 
Let $\Omega_m={g^{-1}\Omega g\cap m^{-1}SWk_g^{-1}}$. One can check that $\nu(\partial\Omega_m)=0$.
To show the claim, we first fix $s\in V$ and $\e>0$. Note that
\begin{eqnarray*}
\Omega_{sm}({a})=\{{k}\in
g^{-1}\Omega g:  {k}k_gB\in {a}^{-1}m^{-1}(s^{-1}V)WB\}.
\end{eqnarray*}
By Lemma \ref{k1}(1), there exists $C>0$ such that
\begin{equation} \label{shahin}
\int_{\Omega_m-\Omega_{sm}({a}) } |f(yg
{k}^{-1}{a}^{-1} (sm)^{-1})| \, d\nu({k}) \le \e
 \quad\text{for all ${a}\in {A}^+_T(C)$.}
\end{equation}
and by Theorem \ref{shah} for any sufficiently large $C>0$ and ${a}\in {A}^+_T(C)$, we have
$$\left| \int_{\Omega_m} f(yg {k}^{-1}{a}^{-1}
(sm)^{-1})  \,d\nu({k}) - \nu(\Omega_m)
\int_{\Lambda\ba L} f\, d\mu\right|<\e.
$$
Hence, \begin{align*} & \left|\int_{{A}^+_T(C)}
\int_{\Omega_{sm}({a}) } f(yg {k}^{-1}{a}^{-1}
(sm)^{-1})\, \xi(\log {a}) d\nu({k}) d{a}\right.\\
& \left.-\int_{A^+_T(C)}\xi(\log a)\,da\cdot\nu(\Omega_m)\cdot \int_{\Lambda\backslash L}f\, d\mu \right|
\\ &\le
\int_{{A}^+_T(C)} \int_{\Omega_m-
\Omega_{sm}({a}) }| f(y g {k}^{-1}{a}^{-1} (sm)^{-1})|
\, \xi(\log {a})
d\nu({k}) d{a}    \\
& +   \int_{{A}^+_T(C)} \left|
 \int_{\Omega_m} f(yg
{k}^{-1}{a}^{-1} (sm)^{-1})\,d\nu(k) - \nu(\Omega_m) \int_{\Lambda\backslash L}f\, d\mu \right|  \, \xi(\log {a})
d{a}
\\ & \le 2\e \int_{{A}^+_T(C)}\xi(\log {a})\ d{a}\le 2\e \int_{{A}^+_T}\xi(\log {a})\ d{a}.
\end{align*}
Thus, by (\ref{eq_ar0}) and (\ref{eq_ar}),
\begin{align*}
&\left|\int_{{A}^+_T}\int_{\Omega_{sm}({a}) } f(yg {k}^{-1}{a}^{-1}
(sm)^{-1})\, \xi(\log {a}) d\nu({k})d{a}\right.\\
&\left.-\rho(B_T(e)) \nu(\Omega_m)
\int_{\Lambda\ba L} f\, d\mu\right|\le 2\varepsilon\rho(B_T(e))+o(\rho(B_T(e)).
\end{align*}
This shows that for every $s\in V$,
\begin{eqnarray*}
\lim_{T\to \infty} \frac{1}{\rho(B_T(e))} \int_{A^+_T}
 \int_{\Omega_{sm}({a}) } f(yg
{k}^{-1}{a}^{-1} (sm)^{-1})\, \xi(\log {a}) d\nu({k}) d{a}\\
 = \nu(\Omega_m)
\int_{\Lambda\backslash L}f\, d\mu .
\end{eqnarray*}
Therefore, (\ref{ddd1}) follows from the Lebesgue dominated convergence theorem.

\vs\vs\noindent\underline{Step 2}: We claim that
for every $m\in M$,
\begin{equation}\label{ddd2}
 \lim_{T\to \infty} \frac{1}{\rho(B_T(e))} \int_{S-\overline{V} } \int_{A^+_T}
\int_{\Omega_{sm}({a}) } f(yg
{{k}}^{-1}{a}^{-1} (sm)^{-1})\, \xi(\log {a}) d\nu({k}) d{a} d\sigma(s)
= 0 .
\end{equation}

Let $s\in S-\overline{V}$ and $\e>0$. By Lemma
\ref{k1}(2), there exists $C>0$ such that
$$\nu(\Omega_{sm}({a}))<\e\quad\text{for all ${a}\in {A}^+_T(C)$.}$$
Hence, by (\ref{eq_ar0}), 
\begin{align*} &\frac{1}{\rho(B_T(e))}  \int_{{A}^+_T(C)}
\int_{\Omega_{sm}({a})} |f(yg{k}^{-1}{a}^{-1}
(sm)^{-1})|\, \xi(\log {a}) d\nu({k}) d{a}\\
& \le  \frac{\e\cdot \sup |f|}{\rho(B_T(e))}\int_{{A}^+_T(C)}\xi(\log {a}) d{a}\le \e\cdot \sup |f|.
\end{align*}
Since $\e>0$ is arbitrary, it follows from (\ref{eq_ar}) that for any
$s\in S -\overline{V} $,
$$\lim_{T\to \infty} \frac{1}{\rho(B_T(e))} \int_{A^+_T}
\int_{\Omega_{sm}({a}) } f(yg
{{k}}^{-1}{a}^{-1} (sm)^{-1})\, \xi(\log {a}) d\nu({k}) d{a} = 0 .$$
Hence, (\ref{ddd2}) follows from the Lebesgue dominated
convergence theorem.

\vs
Since $\sigma(\partial V)=0$, combining (\ref{ddd1}) and (\ref{ddd2}), we deduce that for every $m\in M$,
\begin{align*}
\lim_{T\to \infty} \frac{1}{\rho(B_T(e))} \int_{S}
\int_{A^+_T} \int_{\Omega_{sm}({a}) }
f(yg {{k}}^{-1}{a}^{-1} (sm)^{-1})\, \xi(\log {a}) d\nu({k}) d{a}
d\sigma(s)\\ =\sigma(V)\cdot\nu({g^{-1}\Omega g\cap m^{-1}SWk_g^{-1}})\cdot \int_{\Lambda\backslash L}f\, d\mu .
\end{align*}
Thus, by (\ref{eq_sm}), the Lebesgue dominated
convergence theorem, and Lemma \ref{l_v},
\begin{align*}
&\lim_{T\to \infty} \frac{1}{\rho(B_T(e))} \int_{K}
\int_{A^+_T} \int_{\Omega_{k_1}({a}) }
f(yg {{k}_2}^{-1}{a}^{-1} (k_1)^{-1})\, \xi(\log {a}) d\nu({k}_2) d{a}
d\nu(k_1)\\ &=\sigma(V)\cdot\int_M\nu(g^{-1}\Omega g\cap m^{-1}SWk_g^{-1})\,d\tau(m)\cdot \int_{\Lambda\backslash L}f\, d\mu\\
&=\sigma(V)\nu(g^{-1}\Omega g)\tau(W)\int_{\Lambda\backslash L}f\, d\mu
=\nu(U)\nu_0(\Omega)\int_{\Lambda\backslash L}f\, d\mu.
\end{align*}
Finally, the proposition follows from (\ref{eq_rho_vol}) and (\ref{sum}).
\end{proof}

\vs\begin{proof}[\bf Proof of Theorem \ref{Peq}] 
Let $B$ be a maximal connected split solvable subgroup of $G$ contained in $Q$.
Since $G=KB$, we have $Q=DB$ for $D=K\cap Q$. Then 
\begin{equation}\label{eq_db0}
Q_T(g,\Omega)=DB_T(g,\Omega)
\end{equation}
and for suitable Haar measures $\rho_B$ and $\nu_D$ on $B$ and $D$ respectively.
\begin{equation}\label{eq_db}
\int_{Q} f(q)\,d\rho(q)=\int_D\int_{B} f(db)\,d\rho_B(b)d\nu_D(d),\;\; f\in C_c(Q).
\end{equation}
Hence, Theorem \ref{Peq} for $Q$ follows from Theorem \ref{Peq} for $B$
and the Lebesgue dominated convergence theorem.
Thus, we may assume that $Q=B$.

Let $\e>0$.
One can find a neighborhood $U$ of $e$ in $K$ as in Proposition \ref{weakc} and
functions $f^+, f^-\in C_c(\Lambda\ba L) $ such
that
$$f^-(xk^{-1})\le f(x)\le f^+(xk^{-1})\quad\text{
for all $x\in \Lambda\ba L$ and $k\in U$}$$
 and $$\int_{\Lambda\ba L} |f^+-f^-|\,d\mu\le \e.$$
Put
$$F_T(k)=\frac{1}{\rho(B_T(g))}
\int_{B_T(g, \Omega)} f (y b^{-1}k^{-1})\, d\rho(b),\quad k\in
K.$$
For any $k\in U$,
\begin{multline}\label{ineq}
\frac{1}{\rho(B_T(g))}\int_{B_T(g, \Omega)} f^- (y b^{-1}
k^{-1})\, d\rho(b) \le F_T(e)
 \\ \le
\frac{1}{\rho(B_T(g))}\int_{B_T(g, \Omega)} f^+(y b^{-1}
k^{-1} )\, d\rho(b).
\end{multline}
Integrating over $U\subset K$, we obtain
\begin{multline}\label{conv}
\frac{1}{\rho(B_T(g))} \int_{U} \int_{B_T(g, \Omega)} f^-
(y b^{-1} k^{-1}) \, d\rho(b) d\nu(k)  \le \nu(U) \cdot F_T(e)
\\
\le \frac{1}{\rho(B_T(g))}  \int_{U} \int_{B_T(g, \Omega)}
f^+ (y b^{-1} k^{-1} ) \,d\rho(b)d\nu(k).
\end{multline}
Hence, by Proposition \ref{weakc},
\begin{align*}
&\nu(U)\nu_0(M_0\Omega) \int_{\Lambda\ba L} f^- \, d\mu \le
\nu(U)\cdot \liminf_{T\to\infty} F_T(e) \\
&\le\nu(U)\cdot  \limsup_{T\to\infty} F_T (e) \le
\nu(U)\nu_0(M_0\Omega) \int_{\Lambda\ba L}
f^+\, d\mu.
\end{align*}
This shows that
\begin{align*}
&\nu_0(M_0\Omega)\cdot\left(\int_{\Lambda\ba L} f\, d\mu-\e \right)\le
\liminf_{T\to\infty} F_T(e)\\ &\le \limsup_{T\to\infty} F_T (e)
\le \nu_0(M_0\Omega)\cdot \left(\int_{\Lambda\ba L} f\, d\mu +\e \right),
\end{align*}
and the theorem follows.
\end{proof}

Since it follows from the result of N.~Shah (see \cite[Theorem~1.1]{Sh})
and Ratner's topological rigidity \cite{Ra2} that $\overline{yQ}=\overline{yG_nQ}$ for every $y\in \Lambda\ba L$,
one may expect that Theorem \ref{Peq} holds under the condition $\overline{yQ}=\Lambda\ba L$ as well.

\begin{Lem}\label{l_G}
Let $L$ be a connected semisimple Lie group with finite center and $L=L_cL_n$ the
decomposition of $L$ into the product of compact and noncompact factors.
Suppose that for $y\in\Lambda\ba L$, we have $\overline{yG_n}\supset yL_n$ and $\overline{yG_nQ}=\Lambda\ba L$.
Then the conclusion of Theorem \ref{Peq} holds.
\end{Lem}

\begin{proof}
Let $B\subset Q$ be a maximal connected split solvable subgroup of $G$. Then $Q=DB$ for $D=K\cap Q$.
By Ratner's theorem on orbit closures \cite{Ra2}, the set $\overline{yG_n}$ is a homogeneous space $yG_0$
with a probability $G_0$-invariant measure $\mu_0$, where $G_0$ is a closed connected subgroup of $L$ that contains $L_n$.
Applying Theorem \ref{Peq} to the space $yG_0$ and the subgroup $B$, we deduce that
\begin{equation*}\lim_{T\to \infty} \frac{1}{\rho_B(B_T(g))}\int_{B_T(g,
\Omega)} f(y b^{-1})\, d\rho_B(b)= \nu_0(M_0\Omega)
\int_{yG_0}f \, d\mu_0
\end{equation*}
for every $f\in C_c(\Lambda\ba L)$,
where $\rho_B$ is a right invariant Haar measure on $B$.
Since $\overline{y G_nQ}=\Lambda\ba L$, $L_nQ=L_nD$, and $L_n\subset G_0$, we have 
\begin{equation}\label{eq_1}
\Lambda\ba L=\overline{yG_nD}=yG_0D.
\end{equation}
(Here we used that $D$ is compact.)
By (\ref{eq_db0}), (\ref{eq_db}), and the Lebesgue dominated convergence theorem,
\begin{equation*}\lim_{T\to \infty} \frac{1}{\rho(Q_T(g))}\int_{Q_T(g,
\Omega)} f(y q^{-1})\, d\rho(q)= \nu_0(M_0\Omega)
\int_{\Lambda \ba G}f \, d\tilde\mu
\end{equation*}
for every $f\in C_c(\Lambda\ba L)$, where the measure $\tilde\mu$ is defined by
$$
\int_{\Lambda \ba L}f \, d\tilde\mu=\int_D \int_{yG_0}f(zd^{-1}) \, d\mu_0(z)d\nu_D(d), \quad f\in C_c(\Lambda\ba L).
$$
Since $L_n$ is a normal subgroup of $L$, it is clear that the measure $\tilde\mu$ is invariant under $L_nD$.
This measure corresponds to a Radon measure $\tilde{\mu}$ on $L$ which is right $(L_nD)$-invariant
and left $\hbox{Stab}_L(y)$-invariant. Namely,
$$
\tilde \mu(f)=\int_{\Lambda\ba L} \left(\sum_{\lambda\in\hbox{\small Stab}_L(y)} f(\lambda g)\right) d\tilde \mu(g),\quad f\in C_c(L).
$$
We have a decomposition $\tilde{\mu}=\int_{L_c} \tilde{\mu}_x\, d\omega(x)$
for a Radon measure $\omega$ on $L_c$, where $\tilde\mu_x$ is a right $L_n$-invariant measure on the leaf $xL_n$
for $\omega$-a.e. $x\in L_c$. Since $L_n$ commutes with $x\in L_c$, $\tilde\mu_x$ is left $L_n$-invariant too.
It follows that the measure $\tilde\mu$ is left $L_n$-invariant. Thus, it is left invariant under $\overline{\hbox{Stab}_L(y)L_n}\supset G_0$.
Setting $E=G_0\cap L_c$ and $F=DL_n\cap L_c$, we deduce that
the measure $\omega$ is left $E$-invariant and right $F$-invariant. Note that $EF$ is a closed subset of $L_c$, and
it follows from (\ref{eq_1}) that $yL_nEF=\Lambda\ba L$.
Thus, by the Baire category theorem, the set $EF$ contains an open subset of $L_c$. Since the group $E\times F$ acts transitively
on $EF$, this implies that $EF$ is open in $L_c$. Thus, $L_c=EF$ because $L$ is connected.
By \cite[Theorem~8.32]{Kn}, $\omega$ is a Haar measure on $L_c$.
This implies that $\tilde\mu$ is $L$-invariant, and the lemma follows.
\end{proof}

\section{Distribution of lattice points in sectors and the boundary} \label{pmain}

In this section, we apply Theorem \ref{Peq} in the case when
$\Lambda\ba L=\Gamma\ba G$ in order to deduce Theorems
\ref{thm_main0} and \ref{thm_main}.
We use notations from Theorem \ref{thm_main}.
In particular, $A$ is a split Cartan
subgroup with respect to $K_1$, $g^{-1}K_1g =K_2$, $M_1$ is the
centralizer of $A$ in $K_1$, and $M_2=K_2\cap Q$. 


To simplify notations, we make the following conventions in this section:
\begin{itemize}
\item For $R\subset G$ and $T>0$,
 $$R_T=\{r\in R: d(K_1, K_1 g r) <T\}.$$
\item For $T,C>0$, $\Omega_1\subset K_1$, and $\Omega_2\subset K_2$,
\begin{eqnarray*}
N_T(\Omega_1, \Omega_2)&=& \#(\G_T\cap  g^{-1} K_1 A^+\Omega_1\cap \Omega_2 Q),\\
N^C_T(\Omega_1, \Omega_2)&=& \#(\G_T\cap g^{-1}
K_1 A^+(C)\Omega_1\cap \Omega_2 Q),
\end{eqnarray*}
where $A^+(C)$ is defined as in (\ref{eq_aaa}).
\item For $T,C>0$ and $\Omega\subset K_1$,
\begin{align*} 
Q_T(\Omega)&=Q_T\cap g^{-1}K_1A^+\Omega,\\
Q_T^C(\Omega)&=Q_T\cap g^{-1} K_1 A^+(C) \Omega.
\end{align*}
\end{itemize}
Let $m$ denote the Haar measure on $G$ such that (\ref{cartan}) holds for $K_1$ and $\rho$ the right invariant Haar measure on $Q$ such that
\begin{equation}\label{eq_r}
\int_G f\; dm=\int_{K_2}\int_Q f(kq)\,d\rho(q)d\nu_2(k),\;\;\; f\in C_c(G).
\end{equation}

\begin{Lem}\label{ntc}
For any $C>0$, $\Omega_1\subset K_1$, and $\Omega_2\subset K_2$, 
 $$\lim_{T\to \infty}
\frac{1}{m(G_T)}( N_T(\Omega_1, \Omega_2)- N_T^C(\Omega_1,
\Omega_2)) = 0 .$$
\end{Lem}
\begin{proof}
Fix $D>C>0$  and $\e>0$.
 By Theorem \ref{l_wave}, there exists a neighborhood
  $\Cal O$ of $e$ in $G$ such that
\begin{equation}\label{eq_cal_u}
\Cal O^{-1} g^{-1} K_1A^+(D)K_1\subset g^{-1}K_1A^+(C)K_1.
\end{equation}
In addition, we may choose $\Cal O$ so that
$$
\Gamma\cap \Cal O^{-1} \Cal O
=\{e\}\;\;\;\hbox{and}\;\;\; \Cal O G_T\subset G_{T+\e}\;\;\hbox{for all $T>0$.}
$$
It follows from (\ref{eq_cal_u}) that 
$$ \mathcal{O}\cdot (\Gamma -
g^{-1}K_1A^+(C)K_1) \subset G- g^{-1}
K_1A^+(D)K_1.$$
Thus,
\begin{align*}
&N_T(\Omega_1,\Omega_2)-N^C_T(\Omega_1,\Omega_2)\le \#\{\gamma \in
\Gamma_T - g^{-1}K_1A^+(C)K_1\} \\
&=\frac{1}{m(\Cal O)} m\left(\bigcup_{\gamma\in \Gamma_T - g^{-1}K_1A^+(C)K_1}
 \Cal O\gamma\right)\le \frac{1}{m(\Cal O)}m(G_{T+\e}-g^{-1} K_1A^+(D)K_1)\\
&=\frac{1}{m(\Cal O)}m (K_1(A^+_{T+\e}-A_{T+\e}^+(D))K_1)=o(m(G_{T+\e}))
\end{align*}
by (\ref{eq_B}) and Lemma \ref{cone}. Now the lemma follows from Lemma \ref{vol}.
\end{proof}
The proof of the following lemma is similar and is left to the readers.
\begin{Lem}\label{ntc2}
For any $C>0$ and $\Omega\subset K_1$,
$$\lim_{T\to \infty} \frac{1}{m(G_T)}( \rho(Q_T(\Omega)) -
\rho(Q_T^C(\Omega))) =0 .$$
\end{Lem}


\begin{proof}[\bf Proof of Theorem \ref{thm_main}]
Without loss of generality, $m(\G\ba G)=1$.
It is easy to check that $\nu_1(\partial(M_1\Omega_1))=\nu_2(\partial(\Omega_2M_2))=0$.
Thus, we may assume that $\Omega_1=M_1\Omega_1$ and $\Omega_2=\Omega_2M_2$.

We need to show
that
\begin{equation*}
N_T(\Omega_1, \Omega_2)\sim_{T\to \infty}
\nu_{1}(\Omega_1)\nu_{2}(\Omega_2)m(G_T).
\end{equation*}
Fix any $\e>0$ and $C>0$. Let $U$ be a neighborhood of $e$
in $K_1$ with boundary of measure zero such that $\nu_1(\Omega_1^+- \Omega_1^-)<\e$, where
$$
\Omega_1^+= \bigcup_{u\in U} u\Omega_1 \quad\text{ and }\quad
\Omega_1^-=\bigcap_{u\in U} u^{-1}\Omega_1.
$$
One can check that $\nu_1(\partial\Omega_1^\pm)=0$.

By the strong wavefront lemma (Theorem \ref{l_wave}), there exists a symmetric neighborhood $\mathcal{O}'$ of $e$ in $G$ such that
\begin{equation}\label{o1}
\mathcal{O}'g^{-1}K_1A^+(C)\subset g^{-1}K_1 A^+U\;\;\;\hbox{and}\;\;\; \Cal O' G_T\subset G_{T+\e}\;\;\hbox{for all $T>0$.}
\end{equation}
Set $\mathcal{O}=M_2(\mathcal{O}'\cap Q)M_2$. Note that $\mathcal{O}$ is a symmetric neighborhood of $e$ in $Q$.
Using the fact that $M_2\subset g^{-1}K_1g$, it is easy to check that (\ref{o1}) holds for $\mathcal{O}$ as well.
We may also assume that $\rho(\partial \mathcal{O})=0$.

Let $f$ be the characteristic function of $\Omega_2\Cal O\subset G$.
Since the decomposition $h=h_{K_2}h_Q$ for $h_{K_2}\in K_2$ and $h_Q\in Q$ is uniquely determined
modulo $M_2$ and $M_2\Cal O=\Cal O$, we have
$$
f(h)=\chi_{\Omega_2}(h_{K_2}) \chi_{\Cal O}(h_Q).
$$
We also define a function on $\Gamma\ba G$ by 
$$F (h)=\sum_{\gamma\in \G} f (\gamma h) .$$ 

\vs\vs\noindent\underline{Step 1}: We claim that for any $T>0$,
\begin{equation} N^C_T(\Omega_1,\Omega_2)\le \frac{1}{\rho(\Cal
O)} \int_{Q_{T+\e} (\Omega_1^+)} F(q^{-1})\, d\rho(q), \label{NT1} \end{equation}
We have
\begin{align}\label{last}
\int_{Q_{T + \e} (\Omega_1^{+})}
F(q^{-1}) \, d\rho(q)
&= \int_{Q_{T+\e}(\Omega_1^{+})} \left(\sum_{\gamma\in
\G}\chi_{\Omega_2}(\gamma_{K_2})
 \chi_{\Cal O}(\gamma_Q q^{-1})\right)\, d\rho(q)\\
&=\sum_{\gamma\in \G: \, \gamma_{K_2}\in \Omega_2} \rho(Q_{T+\e}(\Omega_1^{+})\cap \Cal O\gamma_Q).\nonumber
\end{align}
It follows from (\ref{o1})
 that
$$
\Cal OQ_T^C(\Omega_1)\subset Q_{T+\e}(\Omega_1^+).
$$ 
This implies that for every $\gamma\in\G$ such that $\gamma_Q\in
Q^C_T(\Omega_1)$,
$$
\rho(Q_{T+\e}(\Omega_1^{+})\cap \Cal O\gamma_Q)=\rho(\Cal O).
$$
Thus,
\begin{align*}
 \int_{Q_{T+ \e }(\Omega_1^{+})} F(q^{-1})\, d\rho(q) &\ge 
\#\{\gamma\in\Gamma: \gamma_{K_2}\in \Omega_2,\gamma_Q\in Q^C_T(\Omega_1)\}\cdot \rho(\Cal
O)\\
 &=N^C_T(\Omega_1,\Omega_2)\rho(\Cal O).
\end{align*}
This proves (\ref{NT1}).

\vs\vs\noindent\underline{Step 2}: We claim that for any $T>0$,
 \begin{equation}
N_T(\Omega_1, \Omega_2)\ge \frac{1}{\rho(\Cal O)} \int_{Q^C_{T-\e}
(\Omega_1^-)} F(q^{-1})\,
d\rho(q).\label{NT2}
\end{equation} As in (\ref{last}),
\begin{equation}\label{last0}
\int_{Q^C_{T-\e} (\Omega_1^-)}F(q^{-1})
 \, d\rho(q) =\sum_{\gamma\in \G:
\gamma_{K_2}\in \Omega_2} 
\rho(Q^C_{T-\e}(\Omega_1^{-})\cap \Cal O\gamma_Q).
\end{equation}
Since $U\Omega_1^-\subset \Omega_1$, we have by (\ref{o1}),
$$
\Cal O^{-1} Q_{T-\e}^C(\Omega_1^-)\subset Q_{T}(\Omega_1).
$$
Therefore, for $\gamma\in \G$ such that $\gamma_Q\notin
Q_{T}(\Omega_1)$,
$$
\rho(Q^C_{T-\e}(\Omega_1^{-})\cap \Cal O\gamma_Q)=0.
$$
By (\ref{last0}),
\begin{align*}
 \int_{Q^C_{T- \e }(\Omega_1^{-})} F(q^{-1})\, d\rho(q) &\le 
\#\{\gamma\in\Gamma: \gamma_{K_2}\in \Omega_2,\gamma_Q\in Q_T(\Omega_1)\}\cdot \rho(\Cal
O)\\ 
&= N_T(\Omega_1,\Omega_2)\rho(\Cal O).
\end{align*}
This proves (\ref{NT2}).

\vs
Since the boundary of the set $\Omega_2\mathcal{O}$ has measure $0$ (this can be checked using (\ref{eq_r})),
the function $f$ can be approximated by continuous functions with compact support and
Theorem \ref{Peq} can be applied to the function $F$ (see Lemma \ref{l_G}):
\begin{align}\label{eq_peq}
&\lim_{T\to\infty}\frac{1}{\rho(Q_T)}\int_{Q_{T}(\Omega_1^\pm)} F(q^{-1})\,d\rho(q)=\nu_{1}(\Omega_1^\pm)\int_{G/\G} F\,dm\\
&=\nu_{1}(\Omega_1^\pm)\int_{G} f\,dm
=\nu_{1}(\Omega_1^\pm)\nu_2(\Omega_2)\rho(\mathcal{O})\nonumber.
\end{align}
Note that by (\ref{eq_r}), $m(G_T)=\rho(Q_T)$. Combining (\ref{NT2}), Lemma \ref{ntc2}, and (\ref{eq_peq}),
we deduce that
\begin{align*}
&\liminf_{T\to\infty} \frac{N_T(\Omega_1,\Omega_2)}{m(G_T)}\ge
\liminf_{T\to\infty} \frac{1}{m(G_T)\rho(\mathcal{O})}\int_{Q^C_{T-\e}(\Omega_1^-)} F(q^{-1})\,d\rho(q)\\
&\ge\left(\liminf_{T\to\infty} \frac{m(G_{T-\e})}{m (G_T)}\right)\cdot\liminf_{T\to\infty} \frac{1}{\rho(\mathcal{O})\rho(Q_{T-\e})}\int_{Q_{T-\e}(\Omega_1^-)} F(q^{-1})\,d\rho(q)\\
&\ge a(\e)\nu_1(\Omega_1^-)\nu_{2}(\Omega_2)\ge a(\e)(\nu_1(\Omega_1)-\e)\nu_{2}(\Omega_2),
\end{align*}
where $a(\e)$ is defined in Lemma \ref{vol}.
Since $\e>0$ is arbitrary, it follows from Lemma \ref{vol} that
$$
\liminf_{T\to\infty} \frac{N_T(\Omega_1,\Omega_2)}{m(G_T)}\ge
\nu_1(\Omega_1)\nu_2(\Omega_2).
$$
By Lemma \ref{ntc}, (\ref{NT1}), and (\ref{eq_peq}),
\begin{align*}
&\limsup_{T\to\infty}\frac{N_T(\Omega_1,\Omega_2)}{m(G_T)}=\limsup_{T\to\infty}
\frac{N^C_T(\Omega_1,\Omega_2)}{m(G_T)}\\
&\le\limsup_{T\to\infty} \frac{1}{m(G_T)\rho(\mathcal{O})}\int_{Q_{T+\e}(\Omega_1^+)} F(q^{-1})\,d\rho(q)\\
&\le\left(\limsup_{T\to\infty} \frac{m(G_{T+\e})}{m(G_T)}\right)\cdot\limsup_{T\to\infty} \frac{1}{\rho(\mathcal{O})\rho(Q_{T+\e})}\int_{Q_{T+\e}(\Omega_1^+)} F(q^{-1})\,d\rho(q)\\
&\le b(\e)\nu_1(\Omega_1^+)\nu_{2}(\Omega_2)\le b(\e) (\nu_1(\Omega_1)+\e)\nu_{2}(\Omega_2),
\end{align*}
where $b(\e)$ is defined in Lemma \ref{vol}.
Thus, by Lemma \ref{vol},
$$
\limsup_{T\to\infty} \frac{N_T(\Omega_1,\Omega_2)}{m(G_T)}\le
\nu_1(\Omega_1)\nu_2(\Omega_2),
$$
and the theorem is proved.
\end{proof}

\begin{proof}[\bf Proof of Theorem \ref{thm_main0}]
Note that $G$ is a connected semisimple center-free Lie group with no
compact factors, and $K_{x}$ and $K_y$ are maximal compact subgroups of
$G$. Since $G$ acts transitively on $X$ (see, for example, \cite[Theorem IV.3.3]{hel}),
there exists $g\in G$ such that $y=x g$. Then
$K_y=g^{-1} K_x g$. The closed positive Weyl chamber $\Cal W_x$ at
$x$ is of the form $xA^+$, where $A^+$ is a positive Weyl chamber
in a split Cartan subgroup $A$ of $G$ with respect to $K_x$ . The
stabilizer of $b\in X(\infty)$ in $G$ is a parabolic subgroup $Q$ of $G$
(see \cite[Proposition~III.3.8]{gjt}).
In particular, $Q$ contains a maximal connected split solvable subgroup of $G$.
Let $\pi:K_y\to bG$ be the map given by $k\mapsto bk^{-1}$.
Then for $\Omega_1\subset K_x$ and $\Omega_2\subset bG$,
\begin{align*}
&\#\{\gamma\in\Gamma:\, y\gamma\in \mathcal{W}_{x}\Omega_1\cap
\operatorname{B_T}(x),\,
b\gamma^{-1}\in \Omega_2\}\\
&=\#\{\gamma\in\Gamma\cap g^{-1}K_x{A^+}\Omega_1 \cap \pi^{-1}(\Omega_2)Q:\, d(K_x,K_xg\gamma)<T\}.
\end{align*}
Note that $\pi$ maps the probability Haar measure on $K_y$ to the measure $m_{b,y}$.
Using the fact that the map $\pi$ is open, one can check that the set $\pi^{-1}(\Omega_2)$ has boundary of measure zero
if $m_{b,y}(\partial\Omega_2)=0$.
Hence, Theorem \ref{thm_main0} follows from Theorem \ref{main}.
\end{proof}


\section{Distribution of lattice points in bisectors}

Let $G$ be a connected semisimple Lie group with finite center and $G=KA^+K$ a Cartan decomposition of $G$.
To simplify notations, we fix $g\in G$ and
for $\Omega_1,\Omega_2\subset K$ and $T,C>0$, define
\begin{eqnarray*}
G_T&=&\{h\in G: d(K,Kgh)<T\},\\
G_T(\Omega_1,\Omega_2)&=&g^{-1}\Omega_1A^+\Omega_2\cap G_T,\\
N_T(\Omega_1,\Omega_2)&=& \#(\Gamma\cap G_T(\Omega_1,\Omega_2)).
\end{eqnarray*}
If we set $A_T^+=\{a\in A^+:d(K,Ka)<T\}$, then $G_T(\Omega_1,\Omega_2)=g^{-1}\Omega_1 A_T^+\Omega_2$.

Let $G$ be a closed subgroup of a Lie group $L$ and $\Lambda$ a lattice in $L$. Let $m$
be a Haar measure on $G$ such that (\ref{cartan}) holds and $\mu$ the Haar measure on $L$ such that $\mu(\Lambda\ba L)=1$.

\begin{Thm}\label{Geq} 
Suppose that for $y\in \Lambda \ba L$, the orbit $yG_n$ is dense in $\Lambda\ba L$,
where $G_n$ is the product of all noncompact simple factors of $G$.
For any Borel subsets $\Omega_1,\Omega_2\subset K$
with boundaries of measure zero and $f\in C_c(\Lambda\ba L)$,
\begin{equation*}\lim_{T\to \infty} \frac{1}{m(G_T)}
\int_{G_T(\Omega_1,\Omega_2)} f(y h^{-1})\, dm(h)=
\nu(\Omega_1)\nu(\Omega_2)
\int_{\Lambda \ba L}f \, d\mu .
\end{equation*}
\end{Thm}

\begin{proof}
By (\ref{cartan}),
\begin{align}\label{eq_v0}
\int_{G_T(\Omega_1,\Omega_2)} f(y h^{-1})\, dm(h)=\int_{\Omega_1}\int_{A_T^+}\int_{\Omega_2} f(yk_2^{-1}a^{-1}k_1^{-1}g)\,\xi(\log a)dk_2dadk_1\\
=\int_{\Omega_1}\int_{A_T^+}\int_{\Omega_2^{-1}} f(yk_2a^{-1}k_1^{-1}g)\,\xi(\log a)dk_2dadk_1.\nonumber
\end{align}
Since $\nu(\Omega_2^{-1})=\nu(\Omega_2)$, by Theorem \ref{shah}, for every $\e>0$, there exists $C>0$ such that
\begin{equation}\label{eq_v}
\left| \int_{\Omega_2^{-1}} f(y {k}_2{a}^{-1}{k}_1^{-1}g)\,d{k}_2 - \nu(\Omega_2) \int_{\Lambda\ba L} f\, d\mu\right|<\e
\end{equation}
for all ${a}\in {A}^+_T(C)$. Since $m(G_T)=\int_{A_T^+}\xi(\log a)\,da$, it follows from Lemma \ref{cone} that
\begin{equation}\label{eq_vv}
\int_{A_T^+-A^+_T(C)}\xi(\log a)\,da=o(m(G_T))\;\;\;\hbox{as $T\to\infty$.}
\end{equation}
Combining (\ref{eq_v}) and (\ref{eq_vv}), we get
\begin{align*}
&\left| \int_{A_T^+}\int_{\Omega_2^{-1}} f(y {k}_2{a}^{-1}
{k}_1^{-1}g)\,d{k}_2\xi(\log a)da - \nu(\Omega_2)\cdot\int_{A_T^+}\xi(\log a)\,da\cdot
\int_{\Lambda\ba L} f\, d\mu\right|\\
&\le \int_{A^+_T(C)}\left|\int_{\Omega_2^{-1}} f(y {k}_2{a}^{-1}
{k}_1^{-1}g)\,d{k}_2-\nu(\Omega_2)
\int_{\Lambda\ba L} f\, d\mu \right|\xi(\log a)da\\
&+2\sup|f|\int_{A_T^+-A_T^+(C)}\xi(\log a)\,da\le \varepsilon \int_{A^+_T(C)}\xi(\log a)\,da+ o(m(G_T))\\
&\le \e\cdot m(G_T)+ o(m(G_T))\;\;\;\hbox{as $T\to\infty$.}
\end{align*}
This proves that for every $k_1\in K$,
$$
\lim_{T\to\infty}\frac{1}{m(G_T)}\int_{A^+_T}\int_{\Omega_2^{-1}} f(y
{k}_2{a}^{-1}{k}_1^{-1}g)\,\xi(\log a) d{k}_2da =
\nu(\Omega_2)\int_{\Lambda\ba L} f\, d\mu.
$$
Now the statement follows from (\ref{eq_v0}) and the Lebesgue dominated convergence theorem.
\end{proof}

\begin{proof}[\bf Proof of Theorem \ref{th_KK}]
We need to show that $$
N_T(\Omega_1M, M\Omega_2)\sim_{T\to\infty} \nu(\Omega_1M)\nu(M\Omega_2)\cdot\frac{m(G_T)}{m(\G\ba G)}.
$$
We may assume without loss of generality that $\Omega_1M=\Omega_1$, $M\Omega_2=\Omega_2$, and $m(\G\ba G)=1$. 
Since $M$ contains all compact factors of $G$, we may assume that $G$ contains no compact factors.

Fix $\e>0$ and $C>0$.
There exists a neighborhood $U$ of $e$ in $K$ with boundary of measure zero
such that $\nu(\Omega_i^+- \Omega_i^-)<\e$, $i=1,2$, where
\begin{eqnarray*}
\Omega_i^+= \bigcup_{u\in U} \Omega_iu \quad&\text{ and }&\quad
\Omega_i^-=\bigcap_{u\in U} \Omega_iu^{-1}.
\end{eqnarray*}
Note that $\nu(\partial\Omega_i^\pm)=0$, $\Omega_iU\subset \Omega_i^+$, and $\Omega_i^-U\subset \Omega_i$.
By the strong wavefront lemma (Theorem \ref{l_wave}), there exists a neighborhood $\mathcal{O}$ of
$e$ in $G$ such that
\begin{eqnarray}\label{eq_W}
\mathcal{O}^{-1}\Omega_1A^+(C)\Omega_2 &\subset& \Omega_1^+ A^+ \Omega_2^+,\nonumber\\
\mathcal{O}\Omega_1^-A^+(C)\Omega_2^- &\subset& \Omega_1 A^+\Omega_2,\\
\mathcal{O}^{\pm 1}G_T &\subset& G_{T+\e}\quad\text{ for all $T>0$}.
\nonumber
\end{eqnarray}
Let
\begin{eqnarray*}
G_T^C(\Omega_1,\Omega_2)&=&g^{-1}KA^+(C)K\cap G_T(\Omega_1,\Omega_2),\\
N_T^C(\Omega_1,\Omega_2)&=&\#(\G\cap G_T^C(\Omega_1,\Omega_2)).
\end{eqnarray*}
It is not hard to show (see Lemmas \ref{ntc} and \ref{ntc2} for a similar argument) that
\begin{eqnarray}\label{ntc_new}
m(G_T(\Omega_1,\Omega_2)-G_T^C(\Omega_1,\Omega_2))&=&o(m(G_T)),\\
N_T(\Omega_1,\Omega_2)-N^C_T(\Omega_1,\Omega_2)&=&o(m(G_T))\label{ntc_new2}
\end{eqnarray}
as $T\to\infty$.

Let $f\in C_c(G)$ be such that $f\ge 0$, $\hbox{supp}(f)\subset \mathcal{O}$ and $\int_G f\;dm=1$.
Define a function on $\Gamma\ba G$ by
$$F (h)=\sum_{\gamma\in \G} f (\gamma h).$$
Clearly, $\int_{\G\ba G} F\, dm=1$.
We claim that
\begin{eqnarray}
N^C_T(\Omega_1,\Omega_2)\le \int_{G_{T+\e} (\Omega_1^+,\Omega_2^+)}
 F( h^{-1})\, dm(h) \label{NT1_new}
\end{eqnarray}
and
\begin{eqnarray}
N_T(\Omega_1, \Omega_2)\ge \int_{G^C_{T-\e}(\Omega_1^-,\Omega_2^-)}
F( h^{-1}) \,
dm(h).\label{NT2_new}
\end{eqnarray}
First, we observe that
$$
\int_{G_{T + \e} (\Omega_1^{+},\Omega_2^+)} F(h^{-1}) \, dm(h)
=\sum_{\gamma\in \G} \int_{G_{T+\e}(\Omega_1^{+},\Omega_2^+)\gamma^{-1}} f(h^{-1})\,
dm(h).
$$
By (\ref{eq_W}),
$$
\mathcal{O}^{-1}G_T^C(\Omega_1,\Omega_2)\subset G_{T+\e}(\Omega_1^+,\Omega_2^+).
$$
Thus, for every $\gamma\in\G\cap G_T^C(\Omega_1,\Omega_2)$, we have $\mathcal{O}^{-1}\subset G_{T+\e}(\Omega_1^{+},\Omega_2^+)\gamma^{-1}$ and
$$
\int_{G_{T+\e}(\Omega_1^{+},\Omega_2^+)\gamma^{-1}} f(h^{-1})\, dm(h)=1.
$$
This implies (\ref{NT1_new}).

To prove (\ref{NT2_new}), we use that
$$
\int_{G^C_{T - \e} (\Omega_1^{-},\Omega_2^-)} F(h^{-1})\, dm(h)
=\sum_{\gamma\in \G} \int_{G^C_{T-\e}(\Omega_1^{-},\Omega_2^-)\gamma^{-1}} f(h^{-1})\,
dm(h).
$$
By (\ref{eq_W}),
$$
\mathcal{O} G_{T-\e}^C(\Omega_1^-,\Omega_2^-)\subset G_{T}(\Omega_1,\Omega_2).
$$
Therefore, for $\gamma\in \G- G_T(\Omega_1,\Omega_2)$,we have $\mathcal{O}^{-1}\cap G^C_{T-\e}(\Omega_1^{-},\Omega_2^-)\gamma^{-1}=\emptyset$ and
$$
\int_{G^C_{T- \e}(\Omega_1^{-},\Omega_2)\gamma^{-1}} f(h^{-1})\, dm(h)=0.
$$
This proves (\ref{NT2_new}).

By (\ref{NT2_new}), (\ref{ntc_new}), and Theorem \ref{Geq},
\begin{align*}
&\liminf_{T\to\infty} \frac{N_T(\Omega_1,\Omega_2)}{m(G_T)}\ge \liminf_{T\to\infty}\frac{1}{m(G_T)}\int_{G^C_{T-\e}(\Omega_1^-,\Omega_2^-)}F(h^{-1})\,dm(h)\\
&\ge\left(\liminf_{T\to\infty} \frac{m(G_{T-\e})}{m (G_T)}\right)\cdot\liminf_{T\to\infty}\frac{1}{m(G_{T-\e})}\int_{G^C_{T-\e}(\Omega_1^-,\Omega_2^-)} F(h^{-1})\,dm(h) \\
&\ge a(\e)\nu(\Omega_1^-)\nu(\Omega_2^-)\ge a(\e)(\nu(\Omega_1)-\e)(\nu(\Omega_2)-\e),
\end{align*}
where $a(\e)$ is defined in Lemma \ref{vol}.
Since $\e>0$ is arbitrary, it follows from Lemma \ref{vol} that
$$
\liminf_{T\to\infty} \frac{N_T(\Omega_1,\Omega_2)}{m(G_T)}\ge
\nu(\Omega_1)\nu(\Omega_2).
$$
The opposite inequality for $\limsup$ is proved similarly using (\ref{ntc_new2}), (\ref{NT1_new}), and Theorem \ref{Geq}.
\end{proof}

\section{Measure-preserving lattice actions}\label{comme-section}

Let $G$ be a connected semisimple Lie group with finite center and $\Gamma_1$, $\Gamma_2$ lattices in $G$.
Let $L=G\times G$, $\Lambda=\Gamma_1\times \Gamma_2$,
$H=\{(g, g):g\in G\}$, and $Q$ a closed subgroup of $H$ containing a maximal connected
split solvable subgroup of $H$. Fix an invariant Riemannian metric $d$ on the symmetric space $K\ba G$, and for $T>0$ and $g\in G$, define
$$
Q_T(g)=\{(q,q)\in Q:\, d(K,Kgq)<T\}.
$$

\begin{Cor}\label{com} 
Suppose that for $y\in\Lambda\ba L$, the orbit $y Q$ is dense in
$\Lambda\ba L$.
Then for any $f\in C_c(\Lambda\ba L)$,
\begin{equation*}\lim_{T\to \infty} \frac{1}{\rho (Q_T(g))}\int_{Q_T(g)} f(y q^{-1})\, d\rho(q)=
\int_{\Lambda \ba L}f
\, d\mu,
\end{equation*}
where $\rho$ is a right Haar measure on $Q$, and $\mu$ is the $L$-invariant probability measure on $\Lambda\ba L$.
\end{Cor}

\begin{proof}
It suffices to check the conditions of Lemma \ref{l_G}. Namely, we need to show that
\begin{equation}\label{eq_m}
\overline{yH_n}\supset yL_n,
\end{equation}
where $H_n$ and $L_n$ denote the product of all noncompact simple
factors of $H$ and $L$ respectively. We also denote by $H_c$ and $L_c$ the product of all compact simple
factors of $H$ and $L$.
By Ratner's theorem on orbit closures \cite{Ra2}, 
$\overline{yH_n}=yH_0$ for some closed subgroup $H_0$ of $L$ containing $H_n$.
Then $yH_0H_c=\overline{yH}=\Lambda\ba L$, and it follows from the Baire category theorem 
(cf. proof of Lemma \ref{l_G})
that $H_0H_c$ is an open subset of $L$.
Since $H_0H_c$ is also closed, we conclude that $L=H_0H_c$. Now (\ref{eq_m}) follows from the Lemma \ref{l_m} below.
\end{proof}

\begin{Lem}\label{l_m}
Let $S$ be a connected subgroup of $L$ that contains $H_n$. Then $S=(S\cap L_n)(S\cap L_c)$.
\end{Lem}

\begin{proof}
Let $\mathfrak{h}_n\subset \mathfrak{s}\subset \mathfrak{l}=\mathfrak{l}_n\oplus \mathfrak{l}_c$ be the corresponding Lie algebras.
We have decompositions $$\mathfrak{h}_n=\bigoplus_i \mathfrak{h}_n^{i}\quad\hbox{and}\quad \mathfrak{l}_n=\bigoplus_i \mathfrak{l}_n^{i},$$
where $\mathfrak{h}_n^{i}$ and $\mathfrak{l}_n^{i}$ are the simple ideals of $\mathfrak{h}_n$ and $\mathfrak{l}_n$ respectively
such that $\mathfrak{h}_n^{i} \subset \mathfrak{l}_n^{i}$. Note that $\mathfrak{h}_n^{i}$ is a maximal subalgebra of $\mathfrak{l}_n^{i}$.
In particular, $\mathfrak{h}_n^{i}$ is its own normalizer in $\mathfrak{l}_n^{i}$.

It suffices to show that for every $s=(\sum_i s_i )+s_c\in\mathfrak{s}$ with $s_i\in\mathfrak{l}_n^{i}$ and $s_c\in\mathfrak{l}_c$,
we have $s_i,s_c\in\mathfrak{s}$. Clearly,
$$
[\mathfrak{h}_n^{i},s]=[\mathfrak{h}_n^{i},s_i]\subset\mathfrak{s}.
$$
If $[\mathfrak{h}_n^{i},s_i]+\mathfrak{h}_n^{i}\ne \mathfrak{h}_n^{i}$, then $[\mathfrak{h}_n^{i},s_i]+\mathfrak{h}_n^{i}$
generates $\mathfrak{l}_n^{i}$, and hence, $\mathfrak{l}_n^{i}\subset\mathfrak{s}$. Otherwise, $s_i$ normalizes $\mathfrak{h}_n^{i}$,
and it folows that $s_i\in\mathfrak{h}_n^{i}$. In both cases, $s_i\in \mathfrak{s}$. 
Then $s_c=s-\sum_i s_i\in \mathfrak{s}$ too. This proves the lemma.
\end{proof}

Theorem \ref{comme} can be deduced from Corollary \ref{com}, as
it was explained in \cite{GW} or \cite{Oh}.

\section{Application to the Patterson-Sullivan
theory}\label{patterson}

In this section we present the proof of Corollary \ref{Albu},
which is based on Theorem \ref{th_KK} and the standard
Abelian argument. 

Let $G$ be the identity component of the isometry group of the symmetric space $X$
and $K$ the stabilizer of $x$ in $G$. Note that $G$ is a connected semisimple Lie
group with trivial center and with no compact factors, and that $K$ is
a maximal compact subgroup of $G$. We may identify $X$ with $K\ba G$. In particular,
$y=Kg$ for some $g\in G$. Let $A$ be a split Cartan subgroup with
respect to $K$, $\mathfrak{a}$ the Lie algebra of $A$, and $A^+$ a
closed Weyl chamber in $A$.
We use notations (\ref{eq_rt}) and (\ref{eq_aaa}).

By (\ref{eq_em}) and Lemma \ref{p_vol_as},
$$\#\{\gamma\in\Gamma: d(x,y\gamma)< T\}=\#(\Gamma\cap g^{-1}G_T)\sim_{T\to\infty} c_1
T^{(r-1)/2}e^{\delta T}$$
for some $c_1>0$, where $\delta=\max\{2\rho(t):t\in \frak a,\, \|t\|\le 1\}$. 
Applying \cite[Corollary 1b, p.~182]{wid}, we obtain that for some constant
$c_2>0$,
\begin{equation}\label{eq_diverge} \sum_{\gamma\in\Gamma} e^{-s
d(x,y \gamma)}
\sim_{s\to\delta^+}\frac{c_2}{(s-\delta)^{(r+1)/2}}.
\end{equation}
This shows in particular that the constant $\delta$ defined in Lemma
\ref{p_vol_as0} coincides with the critical exponent $\delta_\G$ of
the Dirichlet series (\ref{eq_diverge}).

For $H\in\mathfrak{a}$ and $k\in K$, we define a geodesic ray
$$ \sigma_{H,k}(t)=x\exp(tH)k,\quad t\ge 0,
$$
and denote by $[\sigma_{H,k}]$ the element in $X(\infty)$ corresponding to this ray.


Let $\overline X=X\cup X(\infty)$ be the conic compactification of $X$ (see \cite[Ch.~III]{gjt})
and $\mu_x$ a measure on $\overline X$
which is a limit point of the sequence
$\mu_{x,y,s}$ as $s\to\delta^+$. It will follow from our argument that
$\mu_x$ is unique, and hence, $\mu_{x,y,s}\to\mu_x$ as $s\to\delta^+$.

\begin{Lem}\label{l_ps}
The support of $\mu_x$ is contained in $[\sigma_{v,e}]\cdot
K\subset X(\infty)$, where $v\in\mathfrak{a}$ is the barycenter in $\mathfrak{a}^+$ (see Sec.~\ref{s_vol}).
\end{Lem}
\begin{proof} 
Since the set $y\G\subset X$ is discrete, it follows from (\ref{eq_mxy}) and (\ref{eq_diverge}) that $\mu_x(X)=0$.

Let $\mathcal{S}\subset
A^+$ be an open convex cone centered at the origin that
contains the barycenter of $A^+$ in its interior and $C>0$. Then the set
$$
\mathcal{U}_{\mathcal{S}(C)}=x\mathcal{S}(C)K \cup\{[\sigma_{H,k}]\in
X(\infty) :\, H\in \log(\mathcal{S}(C)),\, k\in K\}
$$
is a neighborhood of $[\sigma_{v,e}]\cdot K$ in $\overline{X}$.
It follows from Lemma \ref{l_hel} that two geodesic rays
$\sigma_{s_1,k_1}$ and $\sigma_{s_2,k_2}$ ($s_1,s_2\in
\mathcal{S}$, $k_1,k_2\in K$) are not equivalent unless $s_1$ and
$s_2$ are collinear. Thus,
$$
[\sigma_{v,e}]\cdot K= \bigcap (\overline{\Cal U}_{{\Cal S}(C)}\cap X(\infty)),
$$
where the intersection is taken for all open convex cones $\Cal S$ in
$A^+$ with center at the origin containing the barycenter and $C>0$.
Hence, it suffices to show that $\mu_x(X)=\mu_x(\overline{X}-\overline{\mathcal{U}}_{\mathcal{S}(C)})=0$.

Fix $\e>0$ and $D> C>0$. For simplicity, write $\Cal U_\Cal S$ for $\Cal U_{\Cal S(C)}$.
Let $\mathcal{T}$ be a cone with the same properties such that $\overline{\mathcal{T}}\subset \mathcal{S}$.
By the strong
wavefront lemma (Theorem \ref{l_wave}), there exists a neighborhood
$\mathcal{O}$ of $e$ in $G$ such that
\begin{align*}
\mathcal{O}^{-1}\cdot g^{-1}K\mathcal{T}(D)K\subset g^{-1}K\mathcal{S}(C)K.
\end{align*}
In addition, we may choose $\mathcal{O}$ such that
$$
\mathcal{O}^{-1}\mathcal{O}\cap \Gamma=\{e\}\;\;\;\hbox{and}\;\;\;\mathcal{O}g^{-1}G_T\subset g^{-1}G_{T+\e}\;\; \hbox{for all $T>0$}.
$$
Let
$$
\beta_\mathcal{S}(T)=\#\{\gamma\in\Gamma: d(x,y\gamma)<T, y\gamma\notin
\overline{\mathcal{U}}_\mathcal{S}\}.
$$
Using (\ref{eq_B}), we have
\begin{align*}
&\beta_\mathcal{S}(T)=\#\left(\G\cap (g^{-1}G_T-g^{-1}K\overline{\mathcal{S}(C)}K)\right)\\
&=\frac{1}{\Vol (\mathcal{O})}\Vol \left(\bigcup\left\{\mathcal{O}\gamma: \gamma\in
 \G\cap (g^{-1}G_T-g^{-1}K\overline{\mathcal{S}(C)}K)\right\}\right)\\
&\ll\Vol \left(\mathcal{O}\cdot
(g^{-1}G_T-g^{-1}K\overline{\mathcal{S}(C)}K)\right)
\le \Vol\left(g^{-1}G_{T+\e}-g^{-1}K\mathcal{T}(D)K\right)\\
&\ll\int_{{A^+_{T+\e}}-\mathcal{T}}
e^{2\rho(\log a)}\, da\ll (T+\e)^{r}e^{\delta' (T+\e)},
\end{align*}
where $\delta'=\max\{2\rho(\log a): a\in
\bar A^+_1-\mathcal{T}\}$. Since $\delta'<\delta$,
it follows from \cite[Theorem 2.1, p.~38]{wid} that the sum $
\int_0^\infty e^{-st}d\beta_\mathcal{S}(t) $ converges for $s=\delta$.
Thus, for $f\in C(\overline{X})$ such that
$\hbox{supp}(f)\subset \overline{X}-\overline{\mathcal{U}}_\mathcal{S}$, 
$$
\sum_{\gamma\in \G} e^{-sd(x,y\gamma)}f(y\gamma)\le \sup
|f| \cdot \int_0^\infty e^{-st}d\beta_\mathcal{S}(t)=O(1)\quad\hbox{as $s\to\delta^+$,}
$$ and $\mu_x(f) =0$ by (\ref{eq_mxy}) and (\ref{eq_diverge}).
This proves that $\mu_x(\overline{X}-\overline{\mathcal{U}}_\mathcal{S})=0$. 
\end{proof}

Now to prove Corollary \ref{Albu}, it suffices to show that for every $f\in C(\overline{X})$ and $k\in K$,
\begin{equation}\label{eq_mx}
\mu_x(f_k)=\mu_x(f),
\end{equation}
where
 $f_k$ is defined as $f_k (z)=f(zk)$ for $z\in \overline{X}$.

Let $C>0$. The set $xA^+(C)K$ is diffeomorphic to the product
$A^+(C)\times (M\backslash K)$, where $M$ is the centralizer of $A$
in $K$, and $\overline{xA^+(C)K}$ is homeomorphic to
$\overline{xA^+(C)}\times (M\backslash K)$.
By Lemma \ref{l_ps}, we
may assume that
$$
\hbox{supp}(f)\subset \overline{xA^+(C)}\times M\backslash K.
$$
Moreover, by the Stone-Weierstrass theorem, it suffices to prove (\ref{eq_mx}) for the functions $f$ of the form
$f(zk)=\phi(z)\psi(k)$ for $\phi\in C(\overline{xA^+(C)})$ and $\psi\in C(M\backslash
K)$.

We claim that \begin{equation}\label{eq_mux}
\mu_x(f)=\phi([\sigma_{v,e}])\int_K\psi\, d\nu,
\end{equation}
where $\nu$ denotes the probability Haar measure on $K$.
Note that this implies (\ref{eq_mx}).
Without loss of generality, $\psi\ge
0$, $\psi\ne 0$.

Let $\varepsilon>0$ and $\mathcal{O}$ be a neighborhood of
$[\sigma_{v,e}]$ in $\overline{xA^+(C)}$ such that
\begin{equation}\label{eq_O}
|\phi(z)-\phi([\sigma_{v,e}])|<\varepsilon\quad \text{ for $z\in
\mathcal{O}$.}
\end{equation}
Let $\chi_i$,
$i=0,\ldots,n$, be a partition of unity on $\overline{xA^+(C)}$ such that
$[\sigma_{v,e}]\notin \hbox{supp}(\chi_i)$ for $i\ne 0$ and
$\hbox{supp}(\chi_0)\subset\mathcal{O}$.
Let $f_i=(\phi\chi_i)\otimes \psi$ for each $0\le i\le n$.
Note that $f=\sum_{i=0}^n f_i$.

For $h\in G$
such that $yh\in xA^+(C)K$,
we denote by $k_h$ the unique element in
$M\backslash K$ such that $yh\in xA^+(C)k_h$. By the properties of Cartan decomposition, the map
$h\mapsto k_h$ is smooth. We have
\begin{eqnarray}\label{eq_pat1}
&&\left|\sum_{\gamma\in\G} e^{-sd(x,y\gamma)}f(y\gamma)-\phi([\sigma_{v,e}])\sum_{\gamma\in \G\cap g^{-1}KA^+(C)K} e^{-sd(x,y\gamma)}\psi(k_\gamma)\right|\\
&\le& \sum_{i=1}^n\sum_{\gamma\in\G} e^{-sd(x,y\gamma)}|f_i(y\gamma)-\phi([\sigma_{v,e}])(\chi_i\otimes\psi)(y\gamma)|\nonumber\\
&+& \sum_{\gamma\in\G\cap g^{-1}KA^+(C)K}
e^{-sd(x,y\gamma)}|f_0(y\gamma)-\phi([\sigma_{v,e}])(\chi_0\otimes\psi)(y\gamma)|.\nonumber
\end{eqnarray}
For each $1\le i\le  n$,
$$\hbox{supp}(|f_i-\phi([\sigma_{v,e}])(\chi_i\otimes\psi)|)\cap
[\sigma_{v,e}]\cdot K=\emptyset.$$ Thus, by Lemma \ref{l_ps},
$$\mu_x(|f_i-\phi(\sigma_{v,e})(\chi_i\otimes\psi)|)=0.$$ It
follows from the definition of $\mu_x$ and (\ref{eq_diverge}) that
as $s\to \delta^+$,
\begin{equation}\label{eq_pat2}
\sum_{\gamma\in\G}
e^{-sd(x,y\gamma)}|f_i(y\gamma)-\phi([\sigma_{v,e}])(\chi_i\otimes\psi)(y\gamma)|=
o((s-\delta)^{-(r+1)/2})
\end{equation}
for each  $1\le i\le  n$.
By (\ref{eq_O}),
$$
|f_0(y \gamma)-\phi([\sigma_{v,e}])(\chi_0\otimes\psi) (y\gamma) |<\varepsilon\cdot
(1\otimes \psi)(y\gamma)=\e\cdot \psi(k_\gamma).
$$
Thus, by Lemma \ref{eq_psi} below,
\begin{equation}\label{eq_pat4}
\sum_{\gamma\in\G}
e^{-sd(x,y\gamma)}|f_0(y\gamma)-\phi([\sigma_{v,e}])(\chi_0\otimes\psi)(y\gamma)|=\varepsilon\cdot
O((s-\delta)^{-(r+1)/2})
\end{equation}
as $s\to\delta^+$.
Combining (\ref{eq_pat1})--(\ref{eq_pat4}), (\ref{eq_diverge}),
and Lemma \ref{eq_psi} below, we deduce that
\begin{eqnarray*}
\limsup_{s\to\delta^+} \mu_{x,y,s}(f)&=&\limsup_{s\to\delta^+}\frac{\sum_{\gamma\in\G} e^{-sd(x,y\gamma)}f(y\gamma)}{\sum_{\gamma\in\G} e^{-sd(x,y\gamma)}}
\\
&\le& (1+\varepsilon\cdot O(1)) \phi([\sigma_{v,e}])\int_K\psi\, d\nu,\\
\liminf_{s\to\delta^+} \mu_{x,y,s}(f)&\ge& (1-\varepsilon\cdot O(1))
\phi([\sigma_{v,e}])\int_K\psi\, d\nu.
\end{eqnarray*}
This proves (\ref{eq_mux}), completing the proof of Corollary
\ref{Albu}, provided that we show the following lemma.

\begin{Lem}\label{eq_psi} For every $\psi\in C(M\ba K)$ and $C>0$,
\begin{equation*} \sum_{\gamma\in \G\cap
g^{-1}KA^+(C)K} e^{-sd(x,y\gamma)}\psi(k_\gamma)\sim_{s\to\delta^+}
\frac{c_2}{(s-\delta)^{(r+1)/2}}\int_K \psi\, d\nu.
\end{equation*}
\end{Lem}
\begin{proof}
 By Theorem \ref{th_KK}, Lemma \ref{ntc}, and Lemma
\ref{p_vol_as}, for
every Borel subset $\Omega\subset M\backslash K$ such that $\nu(\partial\Omega)=0$,
\begin{equation}\label{eq_secc}
\#(\G\cap g^{-1}KA^+_T(C)\Omega)\sim_{T\to\infty} \nu
(\Omega)\cdot \frac{\Vol (G_T)}{\Vol(G/\G)}\sim_{T\to\infty} \nu
(\Omega)\cdot c_1T^{(r-1)/2}e^{\delta T}.
\end{equation}
Since continuous functions on $M\ba K$  can be uniformly approximated by
linear combinations of characteristic functions of Borel subsets
with boundaries of measure zero, we deduce that
\begin{equation*}
\sum_{\gamma\in \Gamma\cap g^{-1}KA_T^+(C)K}
\psi(k_\gamma)\sim_{T\to\infty} \left(\int_K\psi\, d\nu\right)\cdot
c_1T^{(r-1)/2}e^{\delta T}.
\end{equation*}
Thus, the lemma follows from \cite[Corollary 1b, p.~182]{wid}.
\end{proof}

\end{document}